\documentclass[12pt]{amsart}

\usepackage{amssymb}
\usepackage{amsmath}
\usepackage{array}
\usepackage{multirow}
\usepackage{verbatim}
\newcommand{\Lg}{\mbox{$\mathfrak g$}}

\newcommand{\Lh}{\mbox{$\mathfrak h$}}
\newcommand{\Lk}{\mbox{$\mathfrak k$}}
\newcommand{\Lp}{\mbox{$\mathfrak p$}}
\newcommand{\La}{\mbox{$\mathfrak a$}}

\newcommand{\Lm}{\mbox{$\mathfrak m$}}
\newcommand{\Ln}{\mbox{$\mathfrak n$}}

\newcommand{\Lq}{\mbox{$\mathfrak q$}}
\newcommand{\Ls}{\mbox{$\mathfrak s$}}
\newcommand{\Lz}{\mbox{$\mathfrak z$}}

\newcommand{\Ld}{\mbox{$\mathfrak d$}}

\newcommand{\Lf}{\mbox{$\mathfrak f$}}
\newcommand{\Le}{\mbox{$\mathfrak e$}}

\newcommand{\D}{\mbox{$\mathcal D$}}

\newcommand{\Pf}{{\em Proof}. }
\newcommand{\EPf}{\hfill$\square$}

\newcommand{\R}{\mbox{$\mathbb R$}}
\newcommand{\C}{\mbox{$\mathbb C$}}
\newcommand{\Q}{\mbox{$\mathbb H$}}
\newcommand{\K}{\mbox{$\mathbb K$}}

\newcommand{\ad}[1]{\mbox{$\mbox{ad}_{#1}$}}
\newcommand{\Ad}[1]{\mbox{$\mbox{Ad}_{#1}$}}

\newcommand{\im}{\mbox{$\mbox{im}\;$}}
\newcommand{\ann}{\mbox{$\mbox{ann}\;$}}

\newcounter{tab}
\newcommand{\tablab}[1]{\refstepcounter{tab}\label{tab:#1}}

\newcommand{\cref}[1]{Corollary~\ref{#1}}

\newtheorem{thm}{Theorem}[section]
\newtheorem*{thm*}{Theorem}
\newtheorem*{thmmain*}{MAIN THEOREM}
\newtheorem{lem}[thm]{Lemma}

\newtheorem{prop}[thm]{Proposition}
\newtheorem*{prop*}{Proposition}

\theoremstyle{definition}
\newtheorem{defn}{Definition}[section]
\theoremstyle{remark}
\newtheorem{rem}{Remark}[section]
\newtheorem{ex}[rem]{Example}

\def\undertilde#1{\mathord{\vtop{\ialign{##\crcr
$\hfil\displaystyle{#1}\hfil$\crcr\noalign{\kern1.5pt\nointerlineskip}
$\hfil\tilde{}\hfil$\crcr\noalign{\kern1.5pt}}}}}

\title{Semisimple symmetric contact spaces}

\author{Dmitri Alekseevsky}

\address{A. A. Kharkevich Institute for Information Transmission Problems,
  B.Karetnuj per.,19, 127051, Moscow, Russia
  and
University of Hradec Kralove, Faculty of Science,
Rokitanskeho 62, 500 03 Hradec Kralove, Czech Republic} 
\email{dalekseevsky@iitp.ru}

\author{Claudio Gorodski}\thanks{The first author
  acknowledges partial financial support from the Czech Science Foundation
   (grant no.~18-00496S).
  The second author acknowledges partial financial
  support from  CNPq (grant 302882/2017-0) and FAPESP (grant 16/23746-6).
  Part of this work was completed while the authors were visiting
  the University of Florence, for which they would like to thank
  Fabio Podest\`a for his hospitality.}

\address{Instituto de Matem\'atica e Estat\'\i stica, Universidade de 
  S\~ao Paulo, Rua do Mat\~ao, 1010, S\~ao Paulo, SP 05508-090, Brazil}
\email{gorodski@ime.usp.br}

\begin{document}

\begin{abstract}
We classify contact manifolds $(M,\mathcal D)$ which are homogeneous under a
connected semisimple Lie group $G$, and symmetric in the sense
that there exists a contactomorphism of $(M,\mathcal D)$ normalizing $G$,
fixing a point $o$ in $M$
and restricting to minus identity along~$\mathcal D_o$.
\end{abstract}

\maketitle

\vspace{-.2in}

\section{Introduction}

In a short note published in 1990, the first named author
described a simple construction of all contact homogeneous
spaces $M=G/H$ of a given Lie group~$G$ in terms of the orbits of the
coadjoint representation of the group~$G$~\cite{A} analogous to the
Kirillov-Kostant-Souriau construction of symplectic homogeneous spaces.

A \emph{contact structure} on a smooth manifold $M$ of dimension~$2n+1$
is a maximally non-integrable distribution $\mathcal D$ of hyperplanes
in the tangent spaces of $M$; the pair $(M,\mathcal D)$ is then
called a \emph{contact manifold}. Locally $\mathcal D$ can be given as the
field of kernels of a locally defined $1$-form $\theta$, called
a local \emph{contact form}, and the maximal
non-integrability condition refers to the fact that
$\theta\wedge(d\theta)^n$ is never zero. In case $M$ is transversally
orientable,
a contact form can be chosen globally and moreover uniquely up to a
conformal factor.

The automorphism group of a contact structure is infinite-dimensional;
its elements are sometimes called \emph{contactomorphisms}.
A contact manifold $(M,\mathcal D)$ is called \emph{homogeneous}
if it admits a transitive Lie group of contactomorphisms; in this
case the contact structure is called \emph{invariant}.
The problem of describing invariant contact structures on a homogeneous space
of a Lie group can be formulated in terms of Lie algebras.
Fix a Lie group $G$. According to main theorem in~\cite{A}, the invariant contact
structures on a simply-connected homogeneous space $M$ of $G$ fall
into one of the two following disjoint classes:
\begin{enumerate}
\item[1.] $M$ is the universal covering of the projectivization of a
conical coadjoint orbit.
\item[2.] $M$ is the total space of a $1$-dimensional bundle over
a covering of a non-conical coadjoint orbit.
\end{enumerate}
Recall that a nonzero coadjoint orbit of a Lie group $G$ is called
\emph{conical} if together with any point in the orbit also the positive
ray through the point is contained in the orbit. In
subsections~\ref{subsec:conical}
and~\ref{non-conical}
we review in detail and in our context the two
constructions above.

Following~\cite{St,FG,BFG,KZ,Z} we say that a contact manifold $(M,\mathcal D)$
is \emph{symmetric} if each point~$p$ of~$M$ is fixed under an involutive
contactomorphism that restricts to minus identity along~$\mathcal D_p$;
such a contactomorphism is called a \emph{symmetry} at~$p$ (it does
not have to be unique). In the quoted references, the contact manifold
carries additional geometric structure related to $\mathcal D$
and the symmetries are required to preserve it
(e.g. sub-Riemannian, sub-conformal, sub-Hermitian,
Cauchy-Riemann or parabolic structure).
Contrary to the case of Riemannian symmetric spaces,
a  contact   symmetric  space needs not be homogeneous, namely,
the  group generated by   symmetries may act non-transitively on the space.
The  first known example
(an odd-dimensional  projective  space  with  two  deleted points)
was constructed  by Lenka Zalabov\'a~\cite{Z}.
For this reason, in this paper
we further restrict to homogeneous spaces and we
say that a homogeneous contact manifold $(M,\mathcal D)$ of a Lie
group $G$ is  \emph{symmetric} if it admits a symmetry at the
basepoint $o=eH\in M$ that normalizes $G$. Our main result is a
complete classification of simply-connected symmetric homogeneous
contact manifolds of a semisimple Lie group.
For  simplicity,  we call such objects \emph{(semisimple)
symmetric contact spaces}.

Some remarks are in order. Although~\cite{BFG} assumes
there is a compatible Riemannian structure
on the contact distribution, their classification result
includes many homogeneous spaces of non-semisimple Lie groups;
the existence of a Hermitian Cauchy-Riemann structure on the
distribution follows from the compactness of the isotropy group;
their examples associated to simple Lie groups are listed
in Table~\ref{tab:5} below. More generally~\cite{KZ} considers
symmetric Hermitian Cauchy-Riemann structures on
distributions more general than of contact type (of arbitrary codimension),
albeit with no general classification results. The examples
in Tables~\ref{tab:6}, \ref{tab:7} and~\ref{tab:8}
admit an invariant (para-)
Cauchy-Riemann structure that is the
pullback of an invariant (para-) complex structure on the base coadjoint orbit.
The paper~\cite{Z} considers symmetric contact structures associated to
parabolic geometries, which endow the Lie algebra of $G$ with a
so-called contact grading; the flat models of such
parabolic geometries correspond to our examples
of conical type listed in Table~\ref{tab:2}.
Further in~\cite{G,GZ} the authors  study  parabolic contact manifolds
carrying a smooth system of symmetries and give  conditions as to
when such manifolds  are  fibered  over  a reflexion space
in the sense of Loos;
such spaces are related to our examples of non-conical
type listed in Tables~\ref{tab:5}, \ref{tab:6},
\ref{tab:7} and~\ref{tab:8}.
Finally, we believe that the examples in Tables~\ref{tab:3}
and~\ref{tab:1} are too simple and/or
already known, but those in Table~\ref{tab:4} are possibly new.

\subsection{Summary of arguments and results}

Since  we are assuming  the  group   $G$  semisimple,  we may
identify   the   dual  space   $\Lg^*$  of  the Lie    algebra with the Lie
algebra  $\Lg$ via the Killing  form   and   identify coadjoint orbits
with adjoint orbits. Then  the  class~(1)  is identified    with
the projectivization $P\mathrm{Ad}_Ge = \mathrm{Ad}_G (\R e) \subset P \Lg$
of  the  adjoint orbit of a nilpotent element $e \in \Lg$,
and contact manifolds in class~(2)
can then be described in terms of one-dimensional
bundles  over   non-nilpotent
adjoint orbits.

Our  starting point  for  the classification  of  semisimple
symmetric contact spaces in class~(1), namely, those of projective type,
is  the Morozov-Jacobson theorem,
which  allows  to  include a nilpotent  element  $e \in \Lg$ into
a standard basis $(h,e,f)$, called  a $\mathfrak{sl}_2$-triple, of a
$3$-dimensional  subalgebra $\Ls$.
Denote  by  $\Lz=Z_{\mathfrak g}(\Ls)$ the centralizer of $\Ls$,
by $N_{\mathfrak g}(\Ls) = \Ls + \Lz$  the  normalizer
of $\Ls$, and  by  $\Lq$  a reductive complement of $N_{\mathfrak g}(\Ls)$
in $\Lg$. 
The  stability  subalgebra of the projectivized orbit
$M = \mathrm{Ad}_G(\R e)\subset P\Lg$ can be  written  as
\[   \Lh = N_{\mathfrak g}(\R e) = \R h +\R e + \Lz + V \]
where  $V = Z_{\mathfrak q}(e)$ is  the  span  of highest weight vectors in
the $\Ls$-module $\Lq$. Denote  by $W$  the
$\mathrm{ad}_h$-invariant complementary subspace to $V$
such  that  $\Lq = V + W$.
Then $ \Lg = \Lh  + \Lm$ where $\Lm=\R f + W$
is identified with the the tangent  space  $T_oM$ 
and $W$ is identified  with  the  contact hyperplane $\mathcal{D}_o$.
The semisimple element $h$ defines  a gradation $\Lg = \sum_{i\in\mathbb Z} \Lg^i$  where  $\Lg^i$ denotes the  $i$-eigenspace; the largest $m$ such that $\Lg^m\neq0$
is called the \emph{depth} of the gradation.
We also say that  $e$ and the corresponding  $\mathfrak{sl}_2$-subalgebra
$\Ls$ are \emph{odd} (resp.~\emph{even})
if this gradation  is odd  (resp.~even),
which means that there exists (resp.~does not  exist)
an odd number  $j$  with  $\Lg^j \neq 0$.

In  the case in which the Lie  algebra  is  absolutely  simple,
we  prove  the  following  theorem which  describes   all 
symmetric contact spaces  $M = G/H$ 
which  are projectivized orbits  of  odd nilpotent  elements.
A  depth  $2$  gradation
\begin{equation}\label{contact-gradation}
  \Lg = \Lg^{-2} + \Lg^{-1} + \Lg^0 +\Lg^1 + \Lg^2
  \end{equation}
is called  a \emph{contact gradation} if  $\dim \Lg^{\pm 2} =1$
and the bilinear form $\Lg^{-1}\times\Lg^{-1}\to\Lg^{-2}$ induced by the
Lie bracket is non-degenerate.
It is   a real  form  of  the canonical contact gradation
of   the complex Lie algebra  $\Lg^{\mathbb C}$  constructed from 
the  regular   $3$-dimensional  subalgebra
$\Ls^{\mathbb C}(\mu) = \C h_{\mu} + \C e_{\mu} + \C f_{\mu}$ associated  with  a long root $\mu$; this  gradation  is  the  eigenspace decomposition
of  $\mathrm{ad}_{h_{\mu}}$  and $(\Lg^{\mathbb C})^{2}  = \C e_{\mu}$,
$(\Lg^{\mathbb C})^0 = \C h_{\mu} + Z_{\mathfrak g^{\mathbb C}}(\Ls^{\mathbb C}(\mu))$.
All contact gradations on real simple Lie algebras
are known and listed, for instance,
in~\cite[Table~1]{G}.

\begin{thm}\label{main1}
  For an absolutely simple Lie algebra $\Lg$ and
  a contact gradation~(\ref{contact-gradation}),
  the  projectivized orbit $M_e = \mathrm{Ad}_G (\R e)$  of the
  nilpositive  element~$e \in \Lg^{2}$ is a symmetric contact space,
  and these manifolds exhaust all symmetric contact spaces  which
  are projectivized orbits of   odd  nilpotent  elements
  if $\Lg$ is not of   $\sf G_2$-type. For  the normal real
form of $\sf G_2$-type,
  the projectivized orbit of  the nilpositive element of  the  regular
  $3$-subalgebra $\Ls(\beta)$ associated  with  a short  root $\beta$
  is   also a symmetric contact space; the  associated gradation has  depth
  $3$. The complete list  is  given  in Tables~\ref{tab:2} and~\ref{tab:3}.
 \end{thm}

The projectivized orbits  $M_e$  are  called \emph{real adjoint varieties}
and have many  remarkable properties.
The  associated symmetric  decomposition  
has  the  form
\[ \Lg = \Lg_+ + \Lg_-  = (\Ls + \Lz ) + (V + W). \]
In particular, the corresponding manifold $G/G_+$
is a \emph{symmetric para-quaternionic K\"ahler space}~\cite{AC,DJS}. 

In case of the projectivized orbit of a long root vector, the
$\Lz$-modules $V$ and $W$ are isomorphic; 
in Table~\ref{tab:2} we indicate  the  subalgebra
$\Lz$ and the complexification of its representation on $\Lg^{-1}=W$,
which is always of quaternionic type.

Throughout this paper, we indicate the (complex)
representations by their highest weights 
and denote the $i$th fundamental weight
of a complex simple Lie algebra by $\pi_i$
(cf.~\cite[Table~1, p.~224]{GOV}).

\begin{samepage}
{\footnotesize

\[ \begin{array}{|c|c|c|c|c|}
\hline
\mbox{Cartan's type}&\Lg&\Lz& W^{\mathbb C}& \mbox{Depth}  \\
\hline
AI & \mathfrak{sl}_n\R&\mathfrak{gl}_{n-2}\R&\pi_1+\pi_{n-1}&2\\
AIII & \mathfrak{su}_{p,q}&\mathfrak{u}_{p-1,q-1}&\pi_1+\pi_{n-1}&2 \\
BDI & \mathfrak{so}_{p,q}&
\mathfrak{so}_{p-2,q-2}+\mathfrak{sl}_2\R& \pi_1\otimes\pi_1'&2\\
CI & \mathfrak{sp}_{n}\R&\mathfrak{sp}_{n-1}\R&\pi_1&2  \\
DIII & \mathfrak{so}_{2n}^*&
\mathfrak{so}_{2n-4}^*+\mathfrak{so}_3&\pi_1\otimes\pi_1' &2\\
EI & \mathfrak{e}_{6(6)}&\mathfrak{sl}_6\R&\pi_3 &2  \\
EII  &\mathfrak{e}_{6(2)}&\mathfrak{su}_{3,3}&\pi_3 &2\\
EIII  &\mathfrak{e}_{6(-14)}&\mathfrak{su}_{1,5}&\pi_3&2 \\
EV & \mathfrak{e}_{7(7)}&\mathfrak{so}_{6,6}&\pi_6 &2 \\
EVI &\mathfrak{e}_{7(-5)}&\mathfrak{so}_{12}^*&\pi_6 &2\\
EVII & \mathfrak{e}_{7(-25)}&\mathfrak{so}_{2,10}&\pi_6&2\\
EVIII &\mathfrak{e}_{8(8)}&\mathfrak{e}_{7(7)}&\pi_7 &2\\
EIX & \mathfrak{e}_{8(-24)}&\mathfrak{e}_{7(-25)}&\pi_7 &2\\
FI & \mathfrak{f}_{4(4)}&\mathfrak{sp}_3\R&\pi_3  &2\\
G & \mathfrak{g}_{2(2)}&\mathfrak{sl}_2\R&3\pi_1  &2\\
\hline
\end{array}\]
\begin{center}\tablab{2}
  \textsc{Table~\ref{tab:2}: Symmetric contact spaces of an absolutely  simple group $G$ which
  are projectivized orbits of  a highest root vector  (real  adjoint varieties)
  }
\end{center}

}
\end{samepage}

The case of the projectivized orbit of a short root occurs for
$\Lg=\Lg_{2(2)}$ only; the situation is summarized in Table~\ref{tab:3}. 

{\footnotesize

\[\begin{array}{|c|c|c|c|c|}
\hline
\Lg &\Lz & W^{\mathbb C} &  V^{\mathbb C} & \mbox{Depth} \\
\hline
\mathfrak{g}_{2(2)}  & \mathfrak{sl}_2\R & \pi_1+\pi_1+\pi_1& \pi_1 & 3\\
\hline
\end{array}\]
\begin{center}\tablab{3}
\textsc{Table~\ref{tab:3}: Symmetric contact space of  group $G_{2(2)}$ which is   the projectivized orbit of  a short  root vector}
\end{center}

}

The following theorem describes symmetric contact spaces of
which are projectivized orbits of even nilpotent elements
in an absolutely simple Lie algebra $\Lg$.
For an even nilpotent element $e\in\Lg$ and  the corresponding
$\mathfrak{sl}_2\R$-subalgebra $\Ls$, the complexification
$\Ls^{\mathbb C}$ is an $\mathfrak{sl}_2\C$-subalgebra of $\Lg^{\mathbb C}$
whose adjoint action on $\Lg^{\mathbb C}$ admits irreducible components
of odd dimension only, hence it defines an ${\sf SO}_3$-structure
on $\Lg^{\mathbb C}$ in the sense of Vinberg~\cite{V}.
A symmetric contact space which is the projectivized orbit of
an even nilpotent element in $\Lg$ gives rise to an
${\sf SO}_3$-structure on $\Lg^{\mathbb C}$
which is \emph{of symmetric type}, in the sense that the normalizer
$N_{\mathfrak g^{\mathbb C}}(\Ls^{\mathbb C})$ is
a symmetric subalgebra of $\Lg^{\mathbb C}$.
In this case we will prove that the dimensions of the
irreducible components of the adjoint action of $\Ls^{\mathbb C}$
on $\Lg^{\mathbb C}$ do not exceed $5$ and hence the
${\sf SO}_3$-structure is in addition
\emph{short}. In subsection~\ref{sec:even}
we will refer to the classification
of short  ${\sf SO}_3$-structures on complex simple Lie algebras
in~\cite{V} and check which ones are of symmetric type to prove the
following theorem.

\begin{thm}\label{main2}
For an absolutely simple Lie algebra $\Lg$, the symmetric contact 
spaces that are projectivized orbits of even nilpotent elements of $\Lg$
are in bijective correspondence with the isomorphism classes of
short  ${\sf SO}_3$-structures of symmetric type on $\Lg^{\mathbb C}$.
The complete
list in given in Table~\ref{tab:4}.
\end{thm}

{\footnotesize

\[\begin{array}{|c|c|c|c|}
\hline
(\Lg,\Ls + \Lz) & (V+W)^{\mathbb C}  & \mbox{Depth} & \mbox{Conditions}\\
\hline
(\mathfrak{sl}_3\R,\mathfrak{so}_{1,2}) &4\pi_1&4 &-  \\
(\mathfrak{su}_{1,2},\mathfrak{so}_{1,2}) &4\pi_1&4 & - \\
(\mathfrak{so}_{p,q},\mathfrak{so}_{1,2}\oplus\mathfrak{so}_{p-1,q-2})&\pi_1\otimes\pi'_1&2 & p\geq1,\ q\geq2\\
\hline
\end{array}\]
\begin{center}\tablab{4}
  \textsc{Table~\ref{tab:4}: Symmetric contact spaces of an absolutely  simple $G$ which are projectivized orbits of an even nilpotent  element
   }
\end{center}
}

Some remarks about the table are in order.
Recall that $\Ls+\Lz$ embeds into $\Lg$ as a symmetric subalgebra. 
The case $p=q=3$ in Table~\ref{tab:4} gives
\[ (\mathfrak{so}_{3,3},\mathfrak{so}_{1,2}\oplus\mathfrak{so}_{2,1})
=(\mathfrak{sl}_4\R,\mathfrak{sl}_2\R\oplus\mathfrak{sl}_2\R).  \]
We have $\mathfrak{su}^*_4=\mathfrak{sl}_2\Q=\mathfrak{so}_{5,1}$,
$\mathfrak{su}_{2,2}=\mathfrak{so}_{2,4}$,
$\mathfrak{so}_{2,2}=\mathfrak{so}_{1,2}\oplus\mathfrak{so}_{2,1}$,
and
$\mathfrak{so}_4^*=\mathfrak{sk}_2\Q=
\mathfrak{so}_{1,2}\oplus\mathfrak{so}_3$.

The following theorem explains the case of a non-absolutely simple
Lie algebra, in which the possibilities are very limited.

\begin{thm}
For a non-absolutely simple Lie algebra $\Lg$, the only possible
symmetric contact spaces that are projectivized orbits of a nilpotent
element in $\Lg$ are described as follows:
\begin{itemize}
\item[(i)] $\Lg=\mathfrak{sl}_2\R\oplus\mathfrak{sl}_2\R$, $\Lz=0$, 
$\Lg_+=\Ls$ is the diagonal subalgebra and $V+W$ is the skew-diagonal. 
\item[(ii)] $\Lg=\mathfrak{sl}_2\C$, $\Lz=0$, $\Lg_+=\Ls$ is the normal
  real form
of $\Lg$ and $V+W=i\Ls$.
\end{itemize} (cf.~Table~\ref{tab:1}).
\end{thm}

{\footnotesize

  \[\begin{array}{|c|c|c|c|c|}
\hline
\Lg & \Ls & \Lz & (V+W)^{\mathbb C} & \mbox{Depth} \\
\hline
\mathfrak{sl}_2\R\oplus\mathfrak{sl}_2\R & \mbox{Diagonal} & 0 &
2\pi_1 & 2 \\
\mathfrak{sl}_2\C & \mbox{Normal real form} & 0 & 2\pi_1 &2\\
\hline
\end{array}\]
\begin{center}\tablab{1}
  \textsc{Table~\ref{tab:1}: Symmetric contact spaces   $M = G/H$   with  non-absolutely simple Lie groups $G$
  }
\end{center}

}

The last theorem  describes   all  symmetric contact spaces
$M = G/H$ of  a semisimple Lie group $G$
which are associated  with non-nilpotent orbits (non-conical  type).
We prove that $M$ is a canonical contactization of a symplectic
symmetric space~\cite{Bi}.
More precisely, let $(N,\omega)$ be a symplectic manifold.
It is called \emph{quantizable} if there exists
a principal bundle $\pi: P \to M$ with one-dimensional structure
group $A=\R$ or $S^1$ and connection~$\theta$ such that $d\theta=\pi^*\omega$. 
The contact manifold
$(M,\ker\theta)$ is called a \emph{contactization} of $N$,
see~\cite{ACHK}. 

Let $N = \mathrm{Ad}_G\xi = G/K \subset  \Lg$ be
a non-nilpotent adjoint orbit  which is  a  symmetric symplectic manifold
with symmetric  decomposition $\Lg = \Lk + \Lp$ and  symplectic  form
defined  by  the $\mathrm{ad}_{\mathfrak k}$-invariant  closed  2-form
    $\omega(x,y) = d (B \circ\xi) (x,y) := - B(\xi, [x,y])$ for
$x$, $y \in \Lp$, where $B$ is the Killing form of $\Lg$.
Then $\xi$ is semisimple and the
centralizer $K=Z_G(\xi)$ is connected. 
Let $\Lh$ be the $B$-orthogonal complement to~$\xi$
in $\Lk$. Then $\Lh$ is a codimension one ideal of $\Lk$.
Assume that the connected subgroup $H$ generated
by $\Lh$ is a closed subgroup of $K$.
Then the principal $A=K/H$-bundle $M = G/H \to N = G/K$ is a
contactization of $N$. In fact,
the $1$-form $\theta := B \circ\xi$ is  $\mathrm{Ad}_H$ -invariant  and
defines  a contact  form $\theta$ on $M = G/H$ with associated
contact distribution $\mathcal{D}$  which is  the invariant extension
of the  hyperplane $\Lp \subset T_oM$. We prove that every symmetric
contact space $M$ is of this form.

\begin{thm}\label{main4}
  Let $G$ be a connected semisimple Lie group with Lie algebra $\Lg$
  and let $N=\mathrm{Ad}_G\xi=G/K$ be a (semisimple) adjoint orbit
  which is a symplectic symmetric space. Assume that
  the $B$-orthogonal complement to $\xi$ in $\Lk=Z_{\mathfrak g}(\xi)$
  generates a closed subgroup of the centralizer $K=Z_G(\xi)$. 
  Then $N$ admits a
  contactization $M=G/H$ which is a symmetric contact space of $G$.
  Conversely, every (semisimple) symmetric contact space 
  of non-conical type arises in this way.
\end{thm}

The semisimple symplectic symmetric spaces were introduced and
classified by P. Bieliavsky~\cite{Bi}. 
Every simply-connected
symplectic symmetric space of a connected semisimple Lie group
is a direct product of symplectic symmetric spaces of simple Lie groups.
The lists of all possible pairs $(\Lg,\Lk)$ for $\Lg$ simple
are given in Tables~\ref{tab:5}, \ref{tab:6},
\ref{tab:7} and~\ref{tab:8}. The tables are organized and labeled by the
existing type of induced CR structure on the contact
distribution of~$G/H$.

{\footnotesize

\[\begin{array}{|c|c|}
\hline
\Lg & \Lk = \Lh \oplus \R \\
\hline
\mathfrak{su}_n  &\mathfrak{su}_p\oplus\mathfrak{su}_{n-p}\oplus
\mathfrak{u}_1\\
\mathfrak{su}_{p,n-p}  &\mathfrak{su}_p\oplus\mathfrak{su}_{n-p}\oplus
\mathfrak{u}_1\\
\mathfrak{so}^*_{2n} & \mathfrak{su}_n\oplus\mathfrak{so}_2\\
\mathfrak{so}_{2n} & \mathfrak{su}_n\oplus\mathfrak{so}_2\\
\mathfrak{so}_n & \mathfrak{so}_{n-2}\oplus\mathfrak{so}_2\\
\mathfrak{so}_{n-2,2} & \mathfrak{so}_{n-2}\oplus\mathfrak{so}_2\\
\mathfrak{sp}_n(\R) & \mathfrak{su}_n\oplus\mathfrak{so}_2\\
\mathfrak{sp}_n & \mathfrak{su}_n\oplus\mathfrak{so}_2\\
\mathfrak{e}_{6} & \mathfrak{so}_{10}\oplus\mathfrak{so}_2\\
\mathfrak{e}_{6(-14)} & \mathfrak{so}_{10}\oplus\mathfrak{so}_2\\
\mathfrak{e}_7 & \mathfrak{e}_6\oplus\mathfrak{so}_2\\
\mathfrak{e}_{7(-25)} & \mathfrak{e}_6\oplus\mathfrak{so}_2\\
\mathfrak{so}_{p,1} & \mathfrak{so}_{p-1}\oplus\R\\
\hline
\end{array}\]
\begin{center}\tablab5
  \textsc{Table~\ref{tab:5}: Hermitian CR symmetric contact spaces $G/H$}
\end{center}

\[\begin{array}{|c|c|}
\hline
\Lg & \Lk= \Lh \oplus \R  \\
\hline
\mathfrak{su}_{p,q}  &\mathfrak{su}_{r,s}\oplus\mathfrak{su}_{p-r,q-s}\oplus
\mathfrak{so}_2\ \mbox{($r>0$ and $s>0$)}\\
\mathfrak{sl}_{2n}(\R) & \mathfrak{sl}_n(\C)\oplus\mathfrak{so}_2\\
\mathfrak{su}^*_{2n} & \mathfrak{sl}_n(\C)\oplus\mathfrak{so}_2\\
\mathfrak{so}^*_{2n} & \mathfrak{su}_{p,n-p}\oplus\mathfrak{so}_2\
\mbox{($0<p<n$)}\\
\mathfrak{so}^*_{2n} & \mathfrak{so}^*_{2n-2}\oplus\mathfrak{so}_2\\
\mathfrak{so}_{p,q} & \mathfrak{so}_{p-2,q}\oplus\mathfrak{so}_2\ \mbox{($p>2$ and $q>0$)}\\
\mathfrak{so}_{2p,2q} & \mathfrak{su}_{p,q}\oplus\mathfrak{so}_2\ \mbox{($p>0$ and $q>0$)}\\
\mathfrak{sp}_n(\R) & \mathfrak{su}_{p,n-p}\oplus\mathfrak{so}_2\ \mbox{($0<p<n$)}\\
\mathfrak{sp}_{p,q} & \mathfrak{su}_{p,q}\oplus\mathfrak{so}_2\ \mbox{($p>0$ and $q>0$)}\\
\mathfrak{e}_{6(-14)} & \mathfrak{so}_{2,8}\oplus\mathfrak{so}_2\\
\mathfrak{e}_{6(-14)} & \mathfrak{so}^*_{10}\oplus\mathfrak{so}_2\\
\mathfrak{e}_{6(2)} & \mathfrak{so}^*_{10}\oplus\mathfrak{so}_2\\
\mathfrak{e}_{6(2)} & \mathfrak{so}_{4,6}\oplus\mathfrak{so}_2\\
\mathfrak{e}_{7(7)} & \mathfrak{e}_{6(2)}\oplus\mathfrak{so}_2\\
\mathfrak{e}_{7(-5)} & \mathfrak{e}_{6(2)}\oplus\mathfrak{so}_2\\
\mathfrak{e}_{7(-5)} & \mathfrak{e}_{6(-14)}\oplus\mathfrak{so}_2\\
\mathfrak{e}_{7(-25)} & \mathfrak{e}_{6(-14)}\oplus\mathfrak{so}_2\\
\hline
\end{array}\]
\begin{center}\tablab6
  \textsc{Table~\ref{tab:6}: Pseudo-Hermitian CR symmetric contact spaces $G/H$ }
\end{center}

\[\begin{array}{|c|c|}
\hline
\Lg & \Lk= \Lh \oplus \R \\
\hline
\mathfrak{sl}_n(\R) & \mathfrak{sl}_p(\R)\oplus\mathfrak{sl}_{n-p}(\R)\oplus\R\\
\mathfrak{su}^*_{2n} & \mathfrak{su}^*_{2p}\oplus\mathfrak{su}^*_{2n-2p}\oplus\R\\
\mathfrak{su}_{n,n} & \mathfrak{sl}_n(\C)\oplus\R\\
\mathfrak{so}^*_{4n} & \mathfrak{su}^*_{2n}\oplus\R\\
\mathfrak{so}_{p,q} & \mathfrak{so}_{p-2,q}\oplus\R\\
\mathfrak{so}_{n,n} & \mathfrak{sl}_n(\R)\oplus\R\\
\mathfrak{so}_{p,q} & \mathfrak{so}_{p-1,q-1}\oplus\R\\
\mathfrak{sp}_n(\R) & \mathfrak{sl}_n(\R)\oplus\R\\
\mathfrak{sp}_{n,n} & \mathfrak{su}^*_{2n}\oplus\R\\
\mathfrak{e}_{6(6)} & \mathfrak{so}_{5,5}\oplus\R\\
\mathfrak{e}_{6(-26)} & \mathfrak{so}_{1,9}\oplus\R\\
\mathfrak{e}_{7(7)} & \mathfrak{e}_{6(6)}\oplus\R\\
\mathfrak{e}_{7(-25)} & \mathfrak{e}_{6(-26)}\oplus\R\\
\hline
\end{array}\]
\begin{center}\tablab7
  \textsc{Table~\ref{tab:7}:
Para-pseudo-Hermitian CR symmetric contact spaces $G/H$}
\end{center}

\[\begin{array}{|c|c|}
\hline
\Lg & \Lk= \Lh \oplus \R \\
\hline
\mathfrak{sl}_n(\C) & \mathfrak{sl}_p(\C)\oplus\mathfrak{sl}_{n-p}(\C)\oplus\C\\
\mathfrak{so}_{2n}(\C) & \mathfrak{sl}_n(\C)\oplus\C\\
\mathfrak{so}_n(\C) & \mathfrak{so}_{n-2}(\C)\oplus\C\\
\mathfrak{sp}_n(\C) & \mathfrak{sl}_n(\C)\oplus\C\\
\mathfrak{e}_{6}^{\mathbb C}   & \mathfrak{so}_{10}(\C)\oplus\C\\
\mathfrak{e}_7^{\mathbb C} & \mathfrak{e}_6^{\mathbb C}\oplus\C\\
\hline
\end{array}\]
\begin{center}\tablab8
  \textsc{Table~\ref{tab:8}: CR symmetric contact spaces $G/H$ with pseudo-Hermitian and para-pseudo-Hermitian structures}
\end{center}

}

\section{Homogeneous contact manifolds of  a  semisimple Lie  group}\label{contact}

The main result of this work is a classification of symmetric contact
spaces of semisimple Lie groups according to the following definition.

\begin{defn}\label{symm}
  A  homogeneous  contact  manifold  $(M=  G/H,\mathcal D)$ is
  called a \emph{symmetric
contact space} if there is an involutive contactomorphism of
$(M, \mathcal{D})$  that fixes the point $o =eH$,  acts  on  the
contact  subspace $\mathcal{D}_0$ as minus identity, and
normalizes $G$.
\end{defn}

We first give  a     description  of homogeneous   contact manifolds  
$(M = G/H , \mathcal{D})$ of  a  semisimple Lie  group $G$,  
which   follows  from the general  construction  of  homogeneous  
contact manifolds  of  a Lie  group  
given in terms of  coadjoint orbits \cite{A}. 
In the semisimple  case,  we   may
use  the Killing  form to identify
coadjoint  orbits  and adjoint orbits.  

There  are  two  types  of homogeneous contact
manifolds: manifolds of  conical  type (constructed 
as the projectivization  of  a  nilpotent   orbit),  
and manifolds  of  non-conical  type  (constructed as
homogeneous  line  bundles over a non-nilpotent  adjoint orbit).
We   next describe   the   structure of  such manifolds.

\subsection{Manifolds of conical type}\label{subsec:conical} 
We  describe   the invariant  contact structure  on  the    
projectivization  $M = P \mathrm{Ad}_G e  =  G/N_G(\R e)\subset  P\Lg$  of   
a nilpotent orbit
$N =  \mathrm{Ad}_G e $  (that is,   
the orbit of a nilpotent  element ~$e\in\Lg$).
They  are  characterized  as contact manifolds that  have no 
invariant   contact  form. There   are   
only finitely many  such manifolds  for any semisimple Lie group $G$.

By the real version of the Morozov-Jacobson theorem~\cite{CM},
we can find elements $h$, $f\in\Lg$ such that
$(h,e,f)$ is an $\mathfrak{sl}_2$-triple
in $\Lg$, namely, the following bracket relations
are satisfied:
\[ [h,e]=2e,\ [h,f]=-2f,\ [e,f]=h. \]
Denote by $\Ls$ the corresponding $\mathfrak{sl}_2\R$-subalgebra of $\Lg$.
Recall that the finite-dimensional real irreducible representations of
$\mathfrak{sl}_2\R$ are precisely the real forms of complex
irreducible representations of $\mathfrak{sl}_2\C$~\cite{O}.

To  describe   the infinitesimal  structure of  such  homogeneous manifolds,  we   introduce  the    following  notation:
\[ \begin{array}{rl}
\Lz=Z_{\mathfrak g}(\Ls): &\mbox{the centralizer of $\Ls$ in $\Lg$} \\
V: &\parbox[t]{4in}{the span of all highest weight vectors of irreducible
$\Ls$-submodules of $\Lg$ other than $\Ls$} \\
W=\sum_{i>0}\mathrm{ad}_f^iV: &\parbox[t]{4in}{the remaining weight spaces not contained
in~$\Ls$} \\
\Lk=Z_{\mathfrak g}(e)= \mathbb{R}e + \Lz + V: &\mbox{the centralizer of $e$ in $\Lg$} \\
\Lh=N_{\mathfrak g}(\R e)=  \mathbb{R}h + \mathbb{R}e + \Lz + V: &\mbox{the normalizer of the line $\R e$ in $\Lg$} \\
B: &\parbox[t]{4in}{the $\ad{}$-invariant symmetric bilinear form on $\Lg$,
normalized so that $B(e,f)=1$.}
\end{array} \]
Note that $N=G/K$ and $M=G/H$ as homogeneous spaces,
where $H=N_G(\R e)$ is a closed subgroup of $G$ with Lie algebra $\Lh$
and $K=Z_G(e)$ is a closed  codimension one subgroup of $H$  
with Lie algebra~$\Lk$.

The subspace $ \Lm = \R f + W $
is a  complementary subspace to $\Lh$ in $\Lg$  which is identified  with
the tangent space $T_oM  =  \Lg/\Lh$ at  the  point $o = eH$,  
so that we have a (non-reductive) decomposition
\begin{equation} \label{main-decomposition}
\Lg = \Lh + \Lm =  (\mathbb{R}h + \mathbb{R}e + \Lz + V) + (\R f + W).
\end{equation}
\begin{defn} The  decomposition~(\ref{main-decomposition}) is called
the \emph{canonical  decomposition} of~$\Lg$ associated to
the $\mathfrak{sl}_2$-triple~$(h,e,f)$.
\end{defn}

We  denote  by  $\theta = B \circ e$ the linear  form on  $\Lg$
dual  to the vector $e$.
Since $H$ preserves $\theta$ up to a multiple,
the kernel $\ker\theta = \Lh +W$
defines an $H$-invariant   subspace $\mathcal D_0=(\Lh+W)/\Lh\cong W$
of $T_oM$.
We   extend it  to an invariant
codimension one distribution $\mathcal{D}$ in  $M$.
\begin{prop}
 The  distribution $\mathcal{D}$ is  an  invariant   contact  distribution on  the manifold   $M = G/H = P \mathrm{Ad}_G e$ .
\end{prop}
\Pf It is  sufficient  to  check that  $d \theta$ is  not  degenerate on  $W$.  We show  that $\ker d\theta  = \Lk$.
Indeed, if  $x \in \ker d \theta$,  then
$$  0 = d \theta (x, \Lg)= - \theta (\mathrm{ad}_x \Lg) = (\mathrm{ad}^*_{x}\theta)(\Lg)$$
which means  that  $x \in Z_{\mathfrak g}(\theta) = Z_{\mathfrak g}(e) = \Lk$. \EPf

\subsection{Manifolds of non-conical type}\label{non-conical}
We start with a characterization of a conical orbit.
\begin{lem}[\cite{A}]
Let $G$ be a connected Lie group  and let
$\theta\in\Lg^*$ be  a $1$-form.
Denote by~$\Lk$ the centralizer $Z_{\mathfrak g}(\theta)$.
Then the coadjoint orbit $N=\mathrm{Ad}_G^*(\theta)$ is
conical if and only if $\theta(\Lk)=0$.
\end{lem}

\Pf Let $G_{\mathbb R\theta}$ (resp.~$\Lg_{\mathbb R\theta}$)
denote the normalizer of the line
$\R\theta$ in $G$ (resp.~$\Lg$). The action of $G_{\mathbb R\theta}$
(resp.~$\Lg_{\mathbb R\theta}$)
on $\R\theta$ defines a homomorphism $\ell: G_{\mathbb R\theta}\to\R^+$
(resp.~$d \ell:\Lg_{\mathbb R\theta}\to\R$). The surjectivity of $\ell$
is equivalent to that of $ d \ell$, so $N$ is conical if and only
if there exists $z\in\Lg$ such that
$\mathrm{ad}_z^*\theta=\theta$.

On the other hand, the linear map
$L_\theta:X\in\Lg\mapsto \mathrm{ad}_X^*\theta=-\theta\circ\mathrm{ad}_X\in\Lg^*$
is skew-symmetric in the sense that $(L_\theta X)Y=(L_\theta Y)X$ for
all~$X$, $Y\in\Lg$. 
It follows that $\im(L_\theta)=\ann(\ker L_\theta)=\ann(\Lk)$.
This  shows   that conicity of the coadjoint orbit
through~$\theta$ is equivalent to $\theta\in\ann(\Lk)$, so we are
done. \EPf

\medskip

Let $G$ be a connected semisimple Lie group and
fix $\theta\in\Lg^*$ such that $N=\mathrm{Ad}^*_G(\theta)\cong G/K$
is a non-conical orbit, where $K=Z_G(\theta)$. Then $\theta$ does not
vanish identically on $\Lk$, where $\Lk=Z_{\mathfrak g}(\theta)$ is the
Lie algebra of~$K$. Let
\[ \Lh = \ker\theta\cap\Lk. \]
As in Subsection~\ref{subsec:conical}, one sees that $\Lk=\ker d\theta$. It follows
that \[ \theta([\Lk,\Lg])=d\theta(\Lk,\Lg)=0 \]
proving that $\Lh$ is a codimension one ideal of $\Lk$.
Let $H$ be the associated connected
subgroup of $K$. If $H$ is closed, then $M:=G/H$ is a homogeneous space
which is
the total space of a $1$-dimensional bundle over $N$ and which carries a
natural
contact structure:

\begin{prop}   The  $1$-form  $\theta$   defines   an invariant  contact  form  on  the manifold  $M = G/H$.
\end{prop}

\Pf Since $\theta|_{\mathfrak h}\equiv0$,
$\theta$ induces a $\Ad H$-invariant element of
 $\Lg^*/\Lh^*\cong(\Lg/\Lh)^*=T_o^*M$,
so a globally defined invariant $1$-form on $M$.
Note that $d\theta$ is the pull-back
of the Kirillov-Kostant-Souriau form on $N\cong G/K$
under $G/H\to G/K$, so
$\theta$ is a contact form on $M$, whose
associated contact distribution is denoted by $\mathcal D$. \EPf

\medskip

Denote by~$\xi$ the element of $\Lg$ dual to $\theta$
under the Killing form $B$. Note that
$\xi$ can be $B$-isotropic, but we can choose $\eta\in\Lk\setminus\Lh$ such that
$\theta(\eta)=B(\xi,\eta)=1$, $\eta$ generates a closed 
subgroup $C$ of $K$ and $K=C\rtimes H$ (the semi-direct product;
compare~\cite[Thm.~3.1, p.~51]{GOV}). 
We may now write
\begin{equation}\label{decomp-non-con}
 \Lg = \Lk  + \Lp,\quad \Lk = \Lh + \R\eta, 
\end{equation}
where  $\Lp$ is a subspace complementary  to  $\Lh$ in $\ker\theta$.
Note that $d\theta$ is non-degenerate on~$\Lp$.

\section{Symmetric contact manifolds of
conical type}\label{conical}

The main result of this section is Proposition~\ref{involutive-structure}
which reduces  the  problem  of  classification
of  symmetric  contact spaces of  conical  type  of a semisimple
Lie group $G$  to  the description  of
$\mathfrak{sl}_2$-subalgebras  of symmetric type of the Lie  algebra
 $\Lg$.
We  will   use  the notation from subsection~\ref{subsec:conical}
and refer to the canonical decomposition~(\ref{main-decomposition}).

\subsection{Gradation   associated   with an
$\mathfrak{sl}_2$-triple}\label{graded}
The   semisimple element  $h$ of an $\mathfrak{sl}_2$-triple $(h,e,f)$ defines
a gradation of the Lie algebra:
\[ \Lg = \sum_{i\in\mathbb Z} \Lg^i, \]
 where   $\Lg^i$   is  the   $i$-eigenspace  of   $\mathrm{ad}_h$.
The largest index $m$ such that $\Lg^m\neq0$ is called
the \emph{depth} of the gradation. The  gradation  is  called
\emph{odd}  if  $\Lg^i \neq 0$ for  some  odd $i$  and  \emph{even}
otherwise. 

In our setting
\[ \Lg^0=\Lz + \R h + W^0\quad\mbox{where}\quad  W^i=W\cap\Lg^i \]
and
\[ V= \sum_{i>0} V^i\quad\mbox{where}\quad V^i=V\cap \Lg^i. \]

\begin{defn}
  Let $\Lg$ be a semisimple Lie algebra. Consider
  an  $\mathfrak{sl}_2$-triple $(h,e,f)$  of $\Lg$
and  the associated 3-dimensional subalgebra $\Ls$.
  \begin{enumerate}
  \item[(i)] We say that $(h,e,f)$ and $\Ls$ are \emph{odd}
    (resp.~\emph{even}) if the associated $\mathrm{ad}_h$-gradation of~$\Lg$
    is odd (resp.~\emph{even}); 
\item[(ii)] We  say   that  $(h,e,f)$ and $\Ls$ 
  are  \emph{of symmetric  type}    if the normalizer
  $N_{\mathfrak g}(\Ls)$ is a symmetric subalgebra, that is, 
     the reductive  decomposition
 $$   \Lg =   N_{\mathfrak g}(\Ls) + \mathfrak{q} = (\Ls + \Lz) +( V + W)  $$
     is  a  symmetric  decomposition.
     (Note  that  in  this  case   the  homogeneous  manifold
$G/N_G(\Ls)$,    which  is  the   space of   all
$3$-dimensional subalgebras conjugate to
     $\Ls$,     is  a para-quaternionic K\"ahler symmetric  space~\cite{AC}.)
     \end{enumerate}
\end{defn}

The   following  proposition  reduces
the classification  of contact  symmetric  spaces of  conical  type
of a semisimple Lie group $G$  to
the  description  of  $\mathfrak{sl}_2(\R)$-subalgebras    of
symmetric  type of the Lie algebra $\Lg$.

\begin{prop}\label{involutive-structure}
  The universal covering of the
  homogeneous contact  manifold  $M = G/H = \mathrm{Ad}_G (\R e) \subset  P \Lg$  of conical  type  is a  symmetric  contact space
 if  and only if  $e$ can be included  in
an  $\mathfrak{sl}_2$-triple  of  symmetric  type.
\end{prop}

\Pf  Let $M =  G/H = \mathrm{Ad}_G (\R e)$ be   a  symmetric   contact space
which is  the projectivization  of  a nilpotent orbit
and consider the canonical decomposition~(\ref{main-decomposition})
\[ \Lg = \Lh + \Lm =  (\mathbb{R}h + \mathbb{R}e + \Lz + V) + (\R f + W). \]

If $(h,e,f)$ is an $\mathfrak{sl}_2$-triple  of
symmetric  type, then
\[ \Lg = \Lg_+ + \Lg_- = (\Ls + \Lz) + (V + W)   \]
under an involution $s$. Let $\tilde G$ be the
simply-connected Lie group with Lie algebra $\Lg$, and $\tilde H$ the
connected subgroup for $\Lh$; it is known that
$\tilde H= Z_{\tilde G}(e)^0$~\cite[\S~6.1]{CM}. 
Now  $\tilde M=\tilde G/\tilde H$
(almost effective presentation)
is a simply-connected homogeneous contact manifold of conical type
covering $M$. 
The  involutive automorphism $s$ of $\Lg$ integrates to
an involutive automorphism of $\tilde G$
which  preserves the  stability  subgroup $\tilde H$  and  hence
induces an involutive  contactomorphism  of  $\tilde M$
fixing  the point  $o$ and
acting  on $\mathcal{D}_o$ as  $-1$.
Hence $\tilde M$ is  a  symmetric   contact space.

Next  we prove  the converse   statement. Assume that
the universal covering $\tilde M$ of $M$ is a symmetric
contact space. The  symmetry at the basepoint of $\tilde M$ 
induces   an involutive    automorphism $s$ of $\Lg$ with
symmetric  decomposition $\Lg= \Lg_+ + \Lg_-$
such that $s[\Lh]=\Lh$ and~$W\subset\Lg_-$.

 Note that $\Lg^{-2}=\R f + W^{-2}$. Since $d\theta$ is non-degenerate
on $W$, we can find $x\in W^i$, $y\in W^j$ with $i+j=-2$ and
$-\theta([x,y])=d\theta(x,y)=-1$. Now
\[ s[x,y]=s(f+[x,y]_W)=sf-[x,y]_W \]
and
\[ [sx,sy]=[-x,-y]=f+[x,y]_W \]
implying that $sf=f$.

Since the $\pm1$-eigendecomposition
$\Lm=\R f + W$ under $s$ is $\ad h$-invariant, we see
that $\ad h$ and $s$ commute on $\Lm$, so
\[ \ad h\circ s = s\circ \ad h = \ad{sh}\circ s \]
as operators on $\Lm$. Put $z=h-sh\in Z_{\mathfrak h}(f)=\Lz$.
Since $z$ centralizes $W$ and $e$, it also centralizes $V$.
Consider the adjoint action of $\Lz$ on $\Ls+V+W$. The
kernel $\Ln_1$ of this action is an ideal of $\Lg$ contained in $\Lh$,
thus $\Ln_1=0$. We have proved above that $z\in\Ln_1$, hence
$sh=h$.

Now it follows from the theory of $\mathfrak{sl}_2$-triples that
$se=e$~\cite[Prop.~2.1, ch.~6, p.~194]{GOV}. We deduce that
\[ s|_{\mathfrak s}=\mathrm{id},\ \mbox{and thus}\ s(\Lz)=\Lz\ \mbox{and}\
s(V)=V. \]
For $0\neq x\in V$, we have $0\neq [f,x]\in W$ so
\[ -[f,x]=s[f,x]=[sf,sx]=[f,sx] \]
implying that $sx+x$ is an element of $V$ that centralizes~$f$,
hence zero. This proves $s|_V=-1$.

We have already shown that
\[ \Lg_+= \Ls + \Lz_+\quad\mbox{and}\quad \Lg_-= \Lz_-+V + W \]
where $\Lz=\Lz_++\Lz_-$ under~$s$. It only remains to check that $\Lz_-=0$.

We first claim that $\Lz_++\Ls+V+W$ is a subalgebra of $\Lg$ and
\[ \Lg = (\Lz_++\Ls+V+W)+\Lz_- \]
is a reductive decomposition; indeed this follows from
\[ [V+W,V+W]\subset[\Lg_-,\Lg_-]\subset\Lg_+=\Ls+\Lz_+ \]
and
\[ [\Lz_-,V+W]\subset(V+W)\cap \Lg_+ =0. \]
Next we can consider the kernel $\Ln_2$ of the adjoint
representation of $\Lz_++\Ls+V+W$ on $\Lz_-$. Of course
$\Ln_2$ is an ideal of $\Lg$, and we have seen that it
contains $\Ls+V+W$. Since
\begin{align*}
 B(\Lz,\Ls+V+W)&=B(\Lz,\R h+W^0)\\&=B(\Lz, \ad e(\R f+W^{-2}))\\
&=B(\ad e\Lz,\R f+W^{-2})\\&=0,
\end{align*}
the Killing orthogonal $\Ln_2^\perp$ is contained in $\Lz\subset\Lh$
and thus $\Ln_2^\perp=0$. This implies~$\Lz_-=0$, as desired. \EPf

\section{Contact  symmetric  space associated  with   contact gradation  of   a semisimple  Lie  algebra}
\begin{defn}\emph{(\v{C}ap-Slovak)} A   depth  $2$   gradation
  $\Lg =\Lg^{-2} + \Lg^{-1} + \Lg^0 + \Lg^1 + \Lg^2$
  of a real (resp.~complex) semisimple
  Lie algebra $\Lg$ 
 is  called  a  \emph{contact gradation}  (resp.~complex contact gradation) if  
 \begin{enumerate}
   \item[(a)] $\dim \Lg^{-2} =1$; and
   \item[(b)] the skew-symmetric
     bilinear form $\Lg^{-1}\times\Lg^{-1}\to\Lg^{-2}$ induced by the
     Lie bracket is nondegenerate.
     \end{enumerate}
 \end{defn}

It turns out contact gradations can exist only on simple Lie
algebras~\cite[Proposition~3.2.4]{CS}.

A  contact gradation   defines  a  homogeneous manifold   $M = G/G^{\geq 0}$,  where $G$ is    the simply  connected  Lie  group   with  the Lie  algebra $\Lg$
and   $G^{\geq 0} $   the  subgroup  generated  by the  non-negative
subalgebra $\Lg^{\geq  0}  =  \Lg^0 + \Lg^1 + \Lg^2  $.

 \begin{prop}\label{examples-coming-from-contact-gradations}
   The  manifold $M = G/G^{\geq 0}$ associated  with  a  contact   gradation   is  a   symmetric  contact space.
 \end{prop}

  Proof. The   contact   distribution
$\mathcal{D}$  is  defined  as    a natural extension  of  the  isotropy invariant  subspace  $   \Lg^{-1} $ of  the   space $\Lg^{<0} = \Lg^{-2} + \Lg^{-1}$  which is identified  with  the tangent  space  $T_oM = \Lg/\Lg^{\geq 0}$.
   Denote  by  $s$ the  involutive  automorphism     of  $\Lg$   associated  with   symmetric   decomposition
   $\Lg = \Lg^{ev} + \Lg^{odd} = (\Lg^{-2} + \Lg^0 + \Lg^2 )  + (\Lg^{-1} + \Lg^1)     $. It  defines   an  involutive  automorphism  $\sigma$ of
    the    Lie group  $G$ which preserves the   subgroup $G^{\geq0}$.
      Then   the  transformation
      $gG^{\geq 0}  \mapsto  \sigma (g) G^{\geq 0}$   preserves the  contact   distribution
      $\mathcal{D}$, fixes  the  point   $o$ and   acts   as  $-1$  on  $\D_0$,
      hence  defines  a  symmetry.
      $\Box$

      \subsection{Canonical contact gradation of a complex simple Lie algebra
      and adjoint variety}
      A     complex  simple  Lie  algebra $\Lg$     admits  a  canonical  complex
      contact gradation
      associated  with  a highest  root
      (cf.~\cite[Theorem~4.2]{W} or~\cite[Proposition~3.2.4]{CS}).
Let $\Lg$ be  a   complex  simple Lie  algebra     with  a Cartan subalgebra $\La$ and corresponding  root  space  decomposition
    $$   \Lg = \La + \sum_{\alpha \in R} \Lg_{\alpha}.                                       $$
Let $\Pi \subset R$  be a  simple  root   system   and   $R_+$
the associated   system of positive  roots.  The  highest  root $\mu \in  R_+$
defines  a  contact gradation
of  $\Lg$   as   follows. Denote  by $R^0 = \mu^{\perp} \cap R$  (resp., $R^0_+ = R^0 \cap R_+$)   the  roots  (resp., positive  roots) orthogonal  to   $\mu$
with respect to the Killing form, and  set
    $R^1 =  R_+ \setminus (\{ \mu \} \cup R_+^0)$.  For  a   set of  roots  $P \subset R$  we  denote  by   $\Lg(P) = \sum_{\alpha \in P} \Lg_{\alpha}$ the    space which is
    the span  of  the root   spaces  associated  with  roots from~$P$. Then
   $$  \Lg = \Lg_{-\mu}  + \Lg (-R^1) + \Lg(R^0) + \Lg(R^1) + \Lg_{\mu} $$
    is a complex  contact gradation which is  called  the  \emph{canonical
      complex contact gradation}   associated  with  the highest  root. The
    associated (compact) complex 
    homogeneous manifold
    $G/G^{\geq  0} = \mathrm{Ad}_G [\Lg_{\mu}]  $  (here $G$ is the
    simply-connected complex semisimple Lie group with Lie algebra $\Lg$)
    is   the orbit  of  the highest
    weight line  $\Lg_{\mu}$ in  the projectivization  $P\Lg$ of  the Lie  algebra and
    is called    the    \emph{adjoint variety} (it is  the only   closed orbit of  $G$  in  $P \Lg$.)

    \subsection{Contact gradations   of a  real absolutely simple  Lie  algebra}
    Let    $\Lg = \sum_j \Lg^j$ be  a gradation  of  a  complex  semisimple Lie  algebra.  It is    the   eigenspace   decomposition of  the  adjoint operator
    $\mathrm{ad}_h$,      where
    $h \in \Lg$ is the uniquely   defined   element  of $\Lg$ such that
    $\mathrm{ad}_h|_{\mathfrak g^j} = j\cdot  \mathrm{id}$, called   the
    \emph{grading  element}.

    A   real  form   $\Lg^{\sigma}$ of  $\Lg$   defined  by  an anti-involution
    (conjugate-linear involutive  automorphism)  $\sigma$
      is called   \emph{consistent}  with     the gradation  of  $\Lg$  if it
      inherits  a gradation     $\Lg^{\sigma}  = \sum_j (\Lg^\sigma)^j$ or,  equivalently,
      $\sigma(h) =h$. Any gradation  of $\Lg^\sigma$ is
      induced    by    a  gradation  of $\Lg$.

   \subsubsection{The contact  gradations  of  classical  Lie   algebras }\label{examples}
   Herein   we  describe  the canonical  gradations of the classical
   Lie  algebras. At  first, we  consider   the  real and complex Lie
   algebras
   $  \mathfrak{sl}_n(\K) = \mathfrak{sl}(V)$, $\mathfrak{sp}_m (\K) = \mathfrak{sp}(V)$,  where $\K = \R$  or  $\C$  and  $ V = \K^n$, with $n = 2m $ in
   the  symplectic case.
     Consider    the  gradation
     $$V = V^{-1} + V^0 + V^{1}  =  \K p + V^0 + \K  q $$
     where,  in     the   symplectic  case,  $\omega$ is a
     symplectic form, $p$, $q$   are isotropic vectors
     with   $\omega(p,q)= 1$   which  span  a nondegenerate   subspace,
     and  $V^0$ is its orthogonal  complement. It    induces a contact gradation  of  the  Lie  algebras   $\Lg = \mathfrak{sl}(V)$, $\mathfrak{sp}(V)$,
     where
   $\Lg^j = \{ A \in \Lg ,\,  A V^i \subset  V^{i+j}    \}$.
     In matrix  notation  with  respect  to the  decomposition
     $V = V^{-1} + V^0 + V^1$,  the gradation  is  described   as
   $$
    \begin{pmatrix}
    \Lg^0& \Lg^{-1}& \Lg^{-2}\\
    \Lg^1& \Lg^{0}& \Lg^{-1}\\
    \Lg^2& \Lg^{1}& \Lg^{0}\\
    \end{pmatrix}
    $$
   More precisely,    it is  given  by
   \begin{multline}
   \mathfrak{sl}(V) = \K p \otimes q^*  +( p\otimes (V^0)^* + V^0 \otimes q^*)\\ +\mathfrak{s}(\K p\otimes p^* + \mathfrak{gl}(V^0) + \K q \otimes q^*)+
     ( q\otimes (V^0)^* + V^0 \otimes p^*) + \K q \otimes p^* \end{multline}
and
    $$  \mathfrak{sp}(V)= S^2V = \K p^2  +p V^0 +(\K p q + \mathfrak{sp}(V^0))+
      qV^0 + \K q^2 .$$
      We  have identified  $\mathfrak{sl}(V)$ with a codimension one subspace
      of $V\otimes V^*$ (the traceless endomorphisms of $V$), and 
$\mathfrak{sp}(V)$  with  the   space  $S^2V$ of    symmetric  bilinear  forms on  $V^*$  using~$\omega$. We denote 
      by $ab$  the symmetric   product   $\frac{1}{2} (a \otimes b + b \otimes a)$, for $ a$, $b \in V $. 
      In   the  first case,  the  grading   element  is
      $h = -p \otimes p^* +q \otimes q^*$  which   is  included  into
      the $\mathfrak{sl}_2$-triple
      $$   h,\    e =q \otimes p^*,\ f = p \otimes  q^*. $$
 In the  second  case the   triple  is
   $$ h = 2pq,\ e= -q^2,\ f = p^2. $$

 Next  we  fix  in  the   space  $V = \C^n$, where $n >2$,
 a  Hermitian metric  $\gamma$  of  signature  $(k+1, \ell + 1)$ for
 $k+ \ell >0$,   such  that
 the   subspace $V^0$ is  nondegenerate   and orthogonal  to    the
 space   spanned   by  $p$, $q$  with   $\gamma(p,q) =1$.
 Then  the   corresponding    real  form
 $\mathfrak{su}_{k+1, \ell +1}$ is  consistent   with  the
 gradation of  $\mathfrak{sl}(V)$. 
       Relative  to  the    decomposition   $V = \C p + V^0 + \C q$, matrices  from $\mathfrak{su}_{k+1,\ell +1}$ have  the form
        $$
        \begin{pmatrix}
     \lambda + i \mu   & -X^*_- &  i \alpha_- \\
    X_+   &   A - \frac{2i\mu}{k + \ell}       &  X_-   \\
    i \alpha_+& -X_+^* &   -\lambda  + i \mu\\
    \end{pmatrix}
        $$
        where    $   A \in  \mathfrak{su}_{k ,\ell}$,   $\lambda$, $\mu$,
        $\alpha_{\pm} \in  \R$,  $X_{\pm} \in \C^{k ,\ell}$  and
          $X^*_{\pm} = \gamma (X_{\pm}, \cdot )$
        is  the   complex  conjugate covector.

          Now  we describe  the canonical gradation of    complex   orthogonal   and  real    pseudo-orthogonal      Lie  algebras
          $\mathfrak{so}(V)$  where  $V = \C^n$ or  $\R^{k+2,\ell+2}$ with $k+ \ell >0$.  Using  the metric $g$ in  $V$, we identify
           $\mathfrak{so}(V)$     with   the  space  $\Lambda^2V$ of  bivectors.
            The  gradation  of $\mathfrak{so}(V)$ is  induced   by  the gradation $V = V^{-1} + V^0 + V^1$  of  $V$  where  $V^{\pm 1}$  are  isotropic  2-dimensional   subspaces  such that  $V^{-1} + V^1$  is  non degenerate  and  $V^0$        is its orthogonal  complement.
            Then
            $$  \mathfrak{so}(V) = \Lambda^2 V   = \Lambda^2 V^{-1} + V^{-1} \otimes V^0 + (V^{-1} \otimes V^1 + \mathfrak{so}(V^0))+
             V^{0} \otimes V^1  + \Lambda^2 V^1 $$
              is   a  contact gradation .
              We  denote  by  $p$, $p'$  a basis of  $P: =V^{-1}$
              and   by  $q$, $q'$   the  dual basis  of
             $Q = V^1 \simeq (V^{-1})^*$.
Note  that $P \otimes Q \simeq \mathfrak{gl}(P) \simeq \mathfrak{gl}(Q)$   and $\Lg^0 \simeq  \mathfrak{gl}_2(\K)  \oplus \mathfrak{so}(V^0)$.
The grading  element   is  $h =   2(q \wedge p +  q' \wedge  p')$. It is included  into  $\mathfrak{sl}_2$  triple
             $$h , f = 2(p \wedge p'), e = -2(q \wedge  q')$$
which  corresponds to    the   graded subalgebra  $\K p\wedge p' + \K h + \K q \wedge q'$  of   the graded  algebra   $\mathfrak{so} (P + Q) $, isomorphic  to $\mathfrak{so}_4(\C)$   for  $\K = \C$  and  $\mathfrak{so}_{2,2}$  for  $\K = \R$.

Finally, in case $V=\C^n$ with $n=2m$ there is a further real form
$\mathfrak{so}^*_{2m}$ of the Lie
algebra $\mathfrak{so}(V)$ consisting of the elements preserving a
nondegenerate skew-Hermitian form. The gradation of $\mathfrak{so}^*_{2m}$
is induced by the gradation of $V$ as above, and $\mathfrak{so}(P+Q)$
is isomorphic to $\mathfrak{so}^*_4$.

 \subsection{Classification  of  real simple Lie  algebras  which  admit  a contact gradation}

                   Here  we   describe   all  real simple Lie algebras which admits  a  contact  gradation.
                   Note  that  a  contact gradation in a real or
                   complex semisimple Lie algebra $\Lg$ is a 
                   \emph{fundamental gradation}, i.e.~the  negative
                   subalgebra $\Lg^{<0}$  is generated   by  $\Lg^{-1}$. 
                  Any  fundamental gradation of a   complex semisimple  Lie  algebra  $\Lg$ is  associated  with  a  subset $\Pi^1$ of   the
                    system of  simple  roots $\Pi$ and  defined  by  the  condition  that
                     $   \mathrm{deg}\, \Lg_{\alpha} =1$   for  $\alpha \in  \Pi^1$
                      and  $\mathrm{deg}\, \Lg_{\alpha} =0$ for  $\alpha  \in \Pi\setminus \Pi^1$~\cite{AMT}.

                      The  following  result  by Djokovi\'c  (see~\cite{Dj}
                      or~\cite[Prop.~3.8]{AM} or~\cite[\S6.2]{AMT}) 
                      gives  a    description  of all gradations   of a real  form $\Lg^{\sigma}$   of  a
                      complex  semisimple Lie  algebra $\Lg$ in terms of Satake  diagrams.

                      \begin{prop}\label{djokovic-criterion}
                        Let $\Lg$   be  a  complex semisimple Lie  algebra  with a gradation   defined  by  $\Pi^1$
                     and let $\Lg^{\sigma}$ be  a   noncompact  real  form defined  by
                     a Satake  diagram.
                     Then  $\Lg^{\sigma}$  is  consistent  with  the gradation
                     if  and only if    all  nodes in the Satake diagram
                     associated  with  roots  from  $\Pi^1$
                     are  white  and there is  no curved  arrow  which  connect  a   root  from $\Pi^1$  with  a  root which is not  in $ \Pi^1$.
                   \end{prop}

                   The  contact gradation  of $\Lg$   is  defined  by  the   set $\Pi^1$
                   which  consists of the simple  roots  $\alpha_i$  associated  with
                   fundamental  weights    $\pi_i$   which  appear in  the
                   decomposition  of  the  highest root   $\mu $
                   in terms of   fundamental  weights.
                   For  example  for  $A_n$,  $\mu =  \pi_1 +\pi_n $   and
                   $\Pi^1 = \{ \alpha_1, \alpha_{n} \}$.
                   Following~\cite{OV},  we  write   down  the    decomposition of  the highest   root in terms of fundamental  weights for   all complex
                   simple  Lie  algebras in Table~\ref{tab:ov}.
                   
                   \[ \begin{array}{|l|c|}
                     \hline
                     \mathsf{A}_1 &2\pi_1 \\
  \mathsf{A}_n \ \mbox{$(n\geq2)$}  &\pi_1 + \pi_n\\
  \mathsf{B}_n  \ \mbox{$(n\geq2)$}  &\pi_2\\
  \mathsf{C}_n \ \mbox{$(n\geq3)$}  &2 \pi_1\\ 
  \mathsf{D}_n \ \mbox{$(n\geq4)$}  & \pi_2\\
  \mathsf{E}_6  & \pi_6\\
  \mathsf{E}_7  & \pi_6\\
  \mathsf{E}_8  & \pi_1\\
  \mathsf{F}_4  & \pi_4\\
  \mathsf{G}_2  & \pi_2 \\
  \hline
\end{array} \]
\begin{center}\tablab{ov}
  \textsc{Table~\ref{tab:ov}: Highest roots  for  complex  simple Lie  algebras}
  \end{center}

Using Djokovi\'c's criterion (Proposition~\ref{djokovic-criterion})
and  analyzing   the list of Satake  diagrams  of  real absolutely
simple  Lie  algebras, we deduce:

\begin{prop}\label{djokovic-classif}
  The contact gradations of the classical  noncompact
  real absolutely simple  Lie  algebras  are
  exhausted  by  the contact gradations   of
  $\mathfrak{sl}_n\R$, $n >2$,
  $\mathfrak{su}_{k+1,\ell+1}$, $k+\ell >0$,
  $\mathfrak{sp}_n\R$, $n>1$,
  $\mathfrak{so}_{k+2,\ell+2}$, $k + \ell >0$, and
  $\mathfrak{so}^*_{2m}$, $m>2$,
  described in subsection~\ref{examples}. 
  Each exceptional  noncompact real absolutely simple Lie  algebra admits  a
  unique contact gradation, up to conjugation,
  with   the exception of  the Lie  algebras  $EIV$ and
  $FII$, which admit no contact gradations. 
   \end{prop}

\section{Classification   of $\mathfrak{sl}_2(\R)$-subalgebras   of   symmetric  type   in  an absolutely  simple Lie  algebra $\Lg$}

\begin{defn}  An $\mathfrak{sl}_2(\mathbb{R})$-subalgebra
$\Ls$ of  a real semisimple  Lie  algebra $\Lg$
is  called \emph{regular} if its  complexification
$\mathfrak{s}^{\mathbb{C}}$ is   a regular 3-dimensional  subalgebra
$\Ls(\mu)$  associated  with  some root $\mu$ (with respect  to a
 Cartan  subalgebra of $\Lg^{\mathbb{C}}$). Here
$\Ls(\mu)$ is spanned over $\C$ by an
 $\mathfrak{sl}_2$-triple $(h_{\mu},e_\mu,f_\mu)$,
where  $e_{\mu}$, $f_\mu$
are   root vectors.
In this case  may  assume  that $\Ls$ is spanned over $\R$
by $(h_{\mu}, e_{\mu}, e_{-\mu} )$.
  \end{defn}

Assume  that $\Ls $ is  a regular $\mathfrak{sl}_2(\mathbb{R})$-subalgebra
of a real absolutely simple Lie algebra $\Lg$, spanned by
the $\mathfrak{sl}_2$-triple~$(h,e,f)$, such that $\Ls^{\mathbb C}=\Ls(\mu)$,
where $\mu$ is a  long  root of $\Lg^{\mathbb{C}}$ with respect to some
Cartan subalgebra.
Then the associated  gradation of $\Lg$ has the  form
\[  \Lg = \Lg^{-2}  +  \Lg^{-1} + \Lg^{0}  + \Lg^{1} + \Lg^{2},   \]
where  $\Lg^2 = \mathbb{R}e$, $\Lg^{-2} = \mathbb{R}f$ and
$\Lg^0 = \R h + \Lz$, where $\Lz$ is the centralizer of $\Ls$.
This is a contact gradation 
and $\Ls$ is an $\mathfrak{sl}_2(\R)$-subalgebra
of $\Lg$ of symmetric type, with symmetric decomposition
\[\Lg =\Lg^{even}+\Lg^{odd} =(\Ls + \Lz) + (V + W), \]
where $V=\Lg^1$, $W=\Lg^{-1}$. 

Any two roots of the same length
in a complex  simple Lie  algebra  are conjugate
by an inner automorphism, so the classification of regular
$\mathfrak{sl}_2(\R)$-subalgebras of an absolutely simple Lie
algebra $\Lg$ amounts to the description of anti-linear involutions
(real forms) of the complex simple Lie algebra $\Lg^{\mathbb C}$
that preserve a
regular $3$-dimensional subalgebra $\Ls(\mu)$; here $\mu$ is a
fixed root of $\Lg^{\mathbb C}$ in the simply-laced case, but it can be
either a fixed long root or a fixed short root in the multiply-laced
case.

\subsection{Case of odd $\mathfrak{sl}_2(\R)$-subalgebras}\label{sec:odd}

Let $\Ls$ be an odd $\mathfrak{sl}_2(\R)$-subalgebra of symmetric type
of an absolutely simple Lie algebra $\Lg$. We  prove  that  $\Ls$
must be a regular subalgebra. 
By assumption, the Killing
orthogonal decomposition
$\Lg=N_{\mathfrak g}(\Ls)+\Lq$ is a symmetric
decomposition. Moreover, the $\ad h$-gradation $\Lg=\sum_i\Lg^i$ is
odd, that is, $\Lg^j\neq0$ for some odd~$j$.

\begin{lem}\label{symmetric-decomposition-odd-symmetric}
  $N_{\mathfrak g}(\Ls)=  \Ls + \Lz =%
\Lg^{even} = \Lg^{-2} + \Lg^0 + \Lg^2$,  with  $\Lg^{-2} =\R f$, $\Lg^{2} =\R e $,
and $\Lq = \Lg^{odd}$.
\end{lem}

\Pf Owing to the fact that $\Lg^i$ is the $i$-eigenspace of $\mathrm{ad}_h$,
we have
\[ N_{\mathfrak{g}}(\Ls)=\Ls+\Lz\subset\R f +\Lg^0+\R e\subset \Lg^{-2} + \Lg^0 + \Lg^2  \subset \Lg^{even}\]  and
then $0 \neq \Lg^{odd} \subset \Lq$. Since
$[\Lg^{odd}, \Lg^{odd}]  + \Lg^{odd}$   is  a non-trivial ideal of  the
simple  Lie   algebra   $\Lg$, it  coincides  with  $\Lg$.
This  shows  that   $N_{\mathfrak g}(\Ls) = \Lg^{even} = [\Lg^{odd}, \Lg^{odd}]$
and  $\Lg^{odd} = \Lq$. \EPf

\begin{prop}\label{odd-symmetric-is-regular}
  An odd $\mathfrak{sl}_2(\R)$-subalgebra $\Ls$ of symmetric type
must be regular.
  Moreover:
  \begin{enumerate}
    \item if $\Lg$ not of type $\sf G_2$, then
  the complexification $\Ls^{\mathbb C}$ is of the form $\Ls(\mu)$,
  where $\mu$ is a
  long root of $\Lg^{\mathbb C}$ with respect to some Cartan subalgebra,
  and hence the associated gradation of $\Lg$ is a contact
  gradation (of depth~$2$);
  \item if $\Lg$ is of type $\sf G_2$ and $\Ls^{\mathbb C}=\Ls(\mu)$ for a
  short root $\mu$, then the associated gradation of $\Lg$ has
  depth~$3$.
  \end{enumerate}
\end{prop}

\Pf Let $(h,e,f)$ be an $\mathfrak{sl}_2$-triple spanning $\Ls$.
Then $h$ is a semisimple element of $\Lg$ and belongs to a Cartan
subalgebra $\La$, which is necessarily contained in
$Z_{\mathfrak g}(h)=\Lg^0$. Owing to Lemma~\ref{symmetric-decomposition-odd-symmetric},
$\Lg^0=\R h+\Lz$. Now $e$ is a root
vector of $\Lg^{\mathbb C}$ with respect to $\La^{\mathbb C}$, say associated
to the root $\mu$. Then $f$ is a root vector associated to $-\mu$.
This already proves that $\Ls$ is a regular subalgebra. 

Consider
the root decomposition of $\Lg^{\mathbb C}$ with respect to $\La^{\mathbb C}$.
Note that each $(\Lg^i)^{\mathbb C}=(\Lg^{\mathbb C})^i$ for $i\neq0$ is a sum of
root spaces. Fix an ordering of the roots that puts $h$ into the Weyl chamber.
Then $\mu$ is a positive root. Denote the depth of the $\mathrm{ad}_h$-gradation
of~$\Lg$ by~$m$. 
If $\tilde\alpha$ denotes the highest root, then the highest root
space $(\Lg^{\mathbb C})_{\tilde\alpha}\subset(\Lg^m)^{\mathbb C}$.
Recall that $\mathrm{ad}_f^k[(\Lg^{\mathbb C})_{\tilde\alpha}]=(\Lg^{\mathbb C})_{\tilde\alpha-k\mu}$
if $\tilde\alpha-k\mu$ is a root. The main observation now is that
\emph{the length of the $\mu$-chain of roots
through $\tilde\alpha$ can be at most~$4$,
and it equals $4$ if and only if $\Lg^{\mathbb C}$ is of $\sf G_2$-type
and $\mu$ is a short root}~\cite[ch.~VI, \S1, no.3]{bourbaki}.

If the depth $m=2$,
from Lemma~\ref{symmetric-decomposition-odd-symmetric} we see
that 
\[ \Lg =\Lg^{-2}+\Lg^{-1}+\Lg^0+\Lg^1+\Lg^2, \]
where
\[ \Lg^{-2}=\R f,\ \Lg^0=\R h + \Lz,\ \Lg^2=\R e, \]
so $\mu$ is the highest root.

If $m\geq3$, then $V=\Lg^m$ and
$\mathrm{ad}_f^k:(\Lg^m)^{\mathbb C}\to(\Lg^{m-2k})^{\mathbb C}$ is injective
for $k:1,\ldots,m$. By the main
observation above, this implies that $m=3$, $\Lg$ is of $\mathsf{G}_2$-type
and $\mu$ is a short root. \EPf
 
\medskip

\textit{Proof of Theorem~\ref{main1}.}
The projectivized adjoint orbit of a nilpositive element of a
contact gradation is a symmetric contact space due to
Proposition~\ref{examples-coming-from-contact-gradations}. 

Conversely, in view of Proposition~\ref{involutive-structure}
we need to classify odd $\mathfrak{sl}_2(\R)$-subalgebras of $\Lg$
of symmetric type. 

Owing to Proposition~\ref{odd-symmetric-is-regular}, the
contact gradations of real absolutely simple Lie algebras
described in~Proposition~\ref{djokovic-classif} exhaust all the
possibilities, unless we are in case $\sf G_2$.

In the case of the normal real form $\Lg=\Lg_{2(2)}$
we check that the $\mathfrak{sl}_2(\R)$-subalgebra associated to a
short root is of symmetric type and odd. Let $\tilde\alpha$ and
$\beta$ be the highest root and highest short
root of $\Lg^{\mathbb C}$ with respect to $\La^{\mathbb C}$, respectively, where
$\La$ is a Cartan subalgebra of $\Lg$. Then the normalizer
$N_{\mathfrak g}(\Ls(\beta))=\Ls(\beta)
+\Ls(\tilde\alpha)$ coincides with $N_{\mathfrak g}(\Ls(\tilde\alpha))$,
which is already known to be a symmetric subalgebra of $\Lg$.
Hence $\Ls(\beta)$ is of symmetric type.
Moreover $V+W=2\mathrm{Sym}^3(\R^2)$ as
an $\Ls(\beta)$-module, which says that the eigenvalues of
$\mathrm{ad}_h$ are $\pm3$, $\pm1$, that is,
$\Ls(\beta)$ is odd. We obtain Tables~\ref{tab:2}
and~\ref{tab:3}. \EPf

\subsection{Case  of even $\mathfrak{sl}_2(\R)$-subalgebras}\label{sec:even}

Let $\Ls$ be an even $\mathfrak{sl}_2(\R)$-subalgebra
of symmetric type of an absolutely simple Lie algebra $\Lg$.
We  prove   here  that the complexification~$\Ls^{\mathbb C}$ defines a
short ${\sf SO}_3$-structure
on~$\Lg^{\mathbb C}$ in  the  sense of E.~Vinberg~\cite{V}.

By assumption, the Killing orthogonal decomposition
$\Lg=N_{\mathfrak g}(\Ls)+\Lq$ is a symmetric decomposition.
Moreover, the $\ad h$-gradation $\Lg=\sum_i\Lg^i$ is even, that is,
$\Lg^j=0$ for all odd $j$.
Although $\Ls$ need not be regular, we will see that it is not far
from being regular by means of the following concept.

\begin{defn}
  An even $\mathfrak{sl}_2(\R)$-subalgebra $\Ls$ of a real semisimple
  Lie algebra $\Lg$, spanned by an $\mathfrak{sl}_2$-triple
  $(h,e,f)$,
  will be called \emph{short} if the eigenvalues of the endomorphism
  $\mathrm{ad}_h$
  belong to the set $\{0,\pm2,\pm4\}$ or, equivalently, irreducible
  submodules of the $\mathrm{ad}_{\mathfrak s}$-module $\Lg$ have
  dimensions $1$, $3$ or $5$. 
\end{defn}

If $\Ls$ is a short even $\mathfrak{sl}_2(\R)$-subalgebra
of a real semisimple Lie algebra $\Lg$, then the
complexification $\Ls^{\mathbb C}$ clearly is a short
${\sf SO}_3$-structure  on $\Lg^{\mathbb C}$
in the sense of Vinberg~\cite{V}. 

 \begin{prop}\label{even-depth}
 An even $\mathfrak{sl}_2(\R)$-subalgebra $\Ls$ of symmetric type 
 of an absolutely simple Lie algebra $\Lg$ is short.
 \end{prop}

 \Pf
   Let $(h,e,f)$ be a $\mathfrak{sl}_2$-triple spanning $\Ls$. 
  We shall prove that the $\ad h$-gradation
  of $\Lg$ has depth at most $4$.
  
 In fact, the simplicity of $\Lg$ implies $[\Lg_-,\Lg_-]=\Lg_+$. 
  Moreover we may assume that the adjoint action of $\Lg_+$
  on $\Lg_-$ is irreducible. Indeed, otherwise the simplicity
  of $\Lg$ implies that there is an $\mathrm{ad}_{\mathfrak g_+}$-irreducible
  decomposition $\Lg_-=\Lg^{(1)}+\Lg^{(-1)}$ such that
  $\Lg=\Lg^{(-1)}+\Lg^{(0)}+\Lg^{(1)}$ defines a gradation
  of depth $1$, where
  $\Lg^{(0)}=\Lg_+$~\cite[App., Lem.~2 and comments thereafter]{tanaka}.
  Since $\mathfrak{sl}_2\C$ is a factor of $\Lg^{(0)}\otimes\C$,
  the classification of gradations of complex simple Lie
  algebras (e.g.~\cite[p.~297]{CS}) gives that
  $(\Lg\otimes\C,\Lg_+\otimes\C)=(\mathfrak{sl}_{n+1}\C,%
  \mathfrak{s}(\mathfrak{gl}_2\C\oplus\mathfrak{gl}_{n-1}))$,
  but in this case the gradation induced by the semisimple element of
  $\mathfrak{sl}_2\C$ is not even.

  Denote by $\Ls$ the span of~$(h,e,f)$ and consider the
  canonical decomposition~(\ref{main-decomposition}). 
  The adjoint action of $\Lg_+=\Ls\oplus\Lz$ on $\Lg_-$ is irreducible.
  Since $\Lz$ must preserve the $\Ls$-isotypical decomposition of
  $\Lg_-$, there must be only one isotypical component, that is
  $\Lg_-=V+W=P_m+\cdots+ P_m$ as an $\mathrm{ad}_{\mathfrak{s}}$-module,
  where~$P_m=\mathrm{Sym}^m(\R^2)$
  is the real irreducible representation of~$\mathfrak{sl}_2(\R)$
  of dimension~$m+1$.

  Owing to the contact condition, for every $0\neq w_0\in W^0$ there exists
  $w_{-2}\in W^{-2}$ such that $[w_0,w_{-2}]=f\in\Ls\subset\Lg^+$.
If $m\geq6$, then we
reach a contradiction to the Jacobi identity in $\Lg$
as follows. Let $0\neq w_{-4}\in W^{-4}$. 
  The Lie bracket $[w_{-4},[w_0,w_{-2}]]=[w_{-4},f]\in W^{-6}$ is nonzero.
  However, $[w_{-4},w_0]$ and $[w_{-4},w_{-2}]$ are weight vectors
  of weights~$-4$ and~$-6$, respectively, lying in $\Lg^+=\Ls+\Lz$; hence
  they are zero. This contradicts the Jacobi identity applied
  to $w_{-4}$, $w_0$, $w_{-2}$. Therefore $m=2$ or $m=4$, as desired.\EPf

\begin{prop}\label{except}
There are no even $\mathfrak{sl}_2(\R)$-subalgebras of symmetric
type in an absolutely simple Lie algebra $\Lg$ of exceptional type.
\end{prop}

\Pf Let $\Ls$ be an even $\mathfrak{sl}_2(\R)$-subalgebra of symmetric type
of~$\Lg$.
In view of Proposition~\ref{even-depth},
there is an induced short ${\sf SO}_3$-structure on $\Lg^{\mathbb C}$.
We shall see however $N_{\mathfrak g^{\mathbb C}}(\Ls^{\mathbb C})$ can never be a symmetric
subalgebra of $\Lg^{\mathbb C}$, reaching a contradiction. 

Recall that $N_{\mathfrak g^{\mathbb C}}(\Ls^{\mathbb C})=\Ls^{\mathbb C}\oplus
Z_{\mathfrak g^{\mathbb C}}(\Ls^{\mathbb C})$.
We run through the classification
of short structures on exceptional complex simple Lie algebras
given in~\cite[\S2.2]{V} and, in each case,
use~\cite[Table~21]{D} to determine
$\dim Z_{\mathfrak g^{\mathbb C}}(\Ls^{\mathbb C})$. We obtain Table~\ref{tab:9};
the \emph{index}
refers to the Dynkin index of the $3$-dimensional Lie subalgebra as
listed in~\cite{D}.

{\small

\[\begin{array}{|c|c|c|}
\hline
\Lg^{\mathbb C} & \mbox{\sl Index of $\Ls^{\mathbb C}$} &  \dim(\Ls^{\mathbb C}\oplus Z_{\mathfrak g^{\mathbb C}}(\Ls^{\mathbb C}))  \\
\hline
\Lg_2^{\mathbb C} & 4 & 3 \\ \hline
\multirow{2}{*}{$\Lf_4^{\mathbb C}$} & 4 & 11 \\
               & 8 & 17 \\   \hline
\multirow{2}{*}{$\Le_6^{\mathbb C}$} & 4 & 19 \\
               & 8 & 17 \\  \hline
\multirow{4}{*}{$\Le_7^{\mathbb C}$} & 4' & 38 \\
                & 8 & 20 \\
& 3'' & 55 \\
& 7 & 17 \\  \hline
\multirow{2}{*}{$\Le_8^{\mathbb C}$} & 4'' & 81 \\
& 8 & 31 \\
\hline
\end{array}\]
\begin{center}\tablab9
  \textsc{Table~\ref{tab:9}: Short structures on exceptional Lie algebras and their
  centralizers}
\end{center}

}

From the classification of Berger~\cite{B}, we immediately see that
$\dim\Ls^{\mathbb C}\oplus Z_{\mathfrak g^{\mathbb C}}(\Ls^{\mathbb C})$ is not the
dimension of a symmetric subalgebra. \EPf

\medskip

\medskip

\textit{Proof of Theorem~\ref{main2}.}
Owing to 
Propositions~\ref{even-depth}
and~\ref{except}, it suffices to run through the cases of complex simple
Lie algebras $\Lg^{\mathbb C}$ of classical type; from the
list of symmetric subalgebras, select those that contain
an ideal isomorphic to $\mathfrak{sl}_2\C$; and check
which of those induce a short
$\sf SO_3$-structure on $\Lg^{\mathbb C}$.

Recall the classification of short $\sf SO_3$-structures on a
classical Lie algebra $\Lg^{\mathbb C}$~\cite[p.~257]{V}.
In case $\Lg^{\mathbb C}=\mathfrak{sl}_n\C$,
an $\sf SO_3$-structure is determined by an $n$-dimensional representation
$\rho:\mathfrak{sl}_2\C\to\mathfrak{sl}_n\C$, which in turn is characterized
  by the dimensions $n_1$, $n_2,\ldots$ of its irreducible
  components, which must have all the same parity. The $\sf SO_3$-structure
  is short if and only if all the $n_i$'s do not exceed $3$.
  The same holds for $\Lg^{\mathbb C}=\mathfrak{so}_n\C$
  (resp.~$\Lg^{\mathbb C}=\mathfrak{sp}_{n}\C$),
  with the addendum that the number of
  $n_i$'s equal to $2$ (resp.~$3$) must be even.
  In Table~\ref{tab:11}, for each classical complex simple Lie algebra,
  and for each symmetric subalgebra
  containing $\mathfrak{sl}_2(\C)$ as an ideal, 
we list $\rho$ and check Vinberg's criterion
  for a short ${\sf SO}_3$-structure. 

  Table~\ref{tab:11} contains the list of symmetric subalgebras
  of classical complex simple Lie algebras of the form
  $\mathfrak{sl}_2(\C)\oplus\Lz$ (cf.~\cite{B}) and an indication of whether
  the $\mathfrak{sl}_2(\C)$-factor defines a short ${\sf SO}_3$-structure.

  {\footnotesize

\[\begin{array}{|c|c|c|c|c|}
\hline
\Lg^{\mathbb C} & \mbox{Symm subalg}  & \mbox{$\mathfrak{sl}_2\C$ is an ideal} & \rho &
\mbox{Short ${\sf SO}_3$-struct} \\
\hline
\multirow{4}{*}{$\mathfrak{sl}_n\C$}& \multirow{2}{*}{$\mathfrak{so}_n\C$} & \mbox{Only if $n=3$}& n_1=3 & \mbox{Yes}\\ 
                                  &  & \mbox{or $n=4$} & n_1=n_2=2 & \mbox{Yes} \\ \cline{2-5}
&\mathfrak{s}(\mathfrak{gl}_k\oplus\mathfrak{gl}_{n-k}\C) & \mbox{Only if $n\geq3$, $k=2$} &
 \mbox{$n_1=2$, $n_2=\cdots=n_{n-2}=1$} & \mbox{No}\\ \cline{2-5}
& \mathfrak{sp}_n\C\ \mbox{($n\geq4$ even)} & \mbox{No} & - & - \\
 \hline
\multirow{3}{*}{$\mathfrak{so}_n\C$}& \multirow{2}{*}{$\mathfrak{so}_k\C\oplus\mathfrak{so}_{n-k}\C$} & \mbox{Only if $k=3$}&
\mbox{$n_1=3$, $n_2=\cdots=n_{n-3}=1$}& \mbox{Yes}\\ 
                                  &  & \mbox{or $k=4$} & \mbox{$n_1=n_2=2$, $n_3=\cdots=n_{n-4}=1$} & \mbox{No} \\ \cline{2-5}
&\mathfrak{gl}_{n/2}\C\ \mbox{($n\geq6$ even)}& \mbox{No} & - & - \\ \hline
\multirow{2}{*}{$\mathfrak{sp}_{n}\C$} & \mathfrak{sp}_{k}\C\oplus\mathfrak{sp}_{n-k}\C & \mbox{Only if $k=1$} &
    \mbox{$n_1=2$, $n_2=\cdots=n_{2n-2}=1$} & \mbox{No}\\ \cline{2-5}
&\mbox{$\mathfrak{gl}_n\C$ ($n\geq3$)} & \mbox{No} & - & - \\ \hline
\end{array} \]
\begin{center}\tablab{11}
  \textsc{Table~\ref{tab:11}: Classical complex simple Lie algebras
  and their symmetric subalgebras}
\end{center}
}

 Finally, we collect real forms of
$(\mathfrak{sl}_3\C,\mathfrak{so}_3\C)$,
$(\mathfrak{sl}_4\C,\mathfrak{so}_4\C)$ and
$(\mathfrak{so}_n\C,\mathfrak{so}_3\C\oplus\mathfrak{so}_{n-3}\C)$
that give examples and obtain Table~\ref{tab:4}. \EPf

\section{Short $\mathfrak{sl}_2$-subalgebras of symmetric type
in non-absolutely simple, semisimple Lie algebras}
Assume  that   the  Lie  algebra $\Lg$ is not  absolutely  simple, i.e. the  complex Lie algebra  $\Lg^{\mathbb{C}}$
is not  simple.
\begin{prop}\label{simple}
  If $\Lg$ is semisimple but not simple,
  then the  only $\mathfrak{sl}_2$-triple of
  symmetric type is the
  triple $(h+h',e+e,f + f')$  associated  with  the
  diagonal  subalgebra $(\mathfrak{sl}_2\R)^d$
  of  the  Lie  algebra
  $\Lg=\mathfrak{sl}_2\R \oplus\mathfrak{sl}'_2\R$. Here $\Lz=0$,
  $\Lg_+=\Ls=(\mathfrak{sl}_2\R)^d$ and the  associated
  canonical decomposition~(\ref{main-decomposition}) is
  \[ \Lg=(\R(h+h')+\R(e+e')+0+\R(e-e'))+(\R(f+f')+(\R(h-h')+\R(f-f'))). \]
  \end{prop}

\Pf Let $\Ls$ be an $\mathfrak{sl}_2\R$-subalgebra of symmetric type
of $\Lg$. Then $\Lg=\Lg_++\Lg_-$ under an involution~$s$ of $\Lg$,
where $\Lg_+=\Ls+\Lz$ and $\Lz$ is the centralizer of $\Ls$ in $\Lg$.
There is an $s$-invariant decomposition
$\Lg=\Lg_1\oplus\cdots\oplus\Lg_r$
into a direct sum of ideals, where for
each $i=1,\ldots,r$, either the Lie algebra $\Lg_i$
is simple, or $\Lg_i=\Ld\oplus\Ld$ is a sum of two copies
of a simple Lie algebra $\Ld$ and $s(x,y)=(y,x)$ for
$(x,y)\in\Ld\oplus\Ld$. For each $i$, there is an involutive
decomposition $\Lg_i=(\Lg_i)_++(\Lg_i)_-$ under the restriction
of~$s$.

Consider the projection $\pi_i:\Lg\to\Lg_i$ and put
$\Ls_i:=\pi_i(\Ls)$. Since $\Ls$ does not centralize nonzero
elements in $\Lg_-$, $\Ls_i\neq0$ for all $i$ and hence,
by simplicity of $\Ls$, $\pi_i$ defines an isomorphism $\Ls\cong\Ls_i$. Upon
this identification, we can now write
\[ \Ls=\{(x,\ldots,x)\in\Ls_1\oplus\cdots\oplus\Ls_r\;|\;x\in\Ls\}. \]
The centralizer
\[ \Lz=Z_{\mathfrak g}(\Ls)=Z_{\mathfrak{g}_1}(\Ls_1)\oplus\cdots\oplus
Z_{\mathfrak{g}_r}(\Ls_r) \]
and $Z_{\mathfrak{g}_i}(\Ls_i)\cap \Ls_i=\{0\}$ for all $i$. Since
$\Ls$ is a proper subset of $\Ls_1\oplus\cdots\oplus\Ls_r$ in case $r\geq2$,
the condition $\Lg_+=\Ls\oplus\Lz$ forces $r=1$.

Since $\Lg=\Lg_1$ is assumed non-simple,
$\Lg=\Ld\oplus\Ld$ where $\Ld\cong\Lg_+$ is simple.
Since $\Lg_+=\Ls\oplus\Lz$ is a sum of ideals,
this implies $\Lz=0$ and $\Ld=\Ls\cong\mathfrak{sl}_2\R$, as we wished. \EPf

\begin{prop}\label{abs-simple}
  If $\Lg$ is a complex simple Lie algebra viewed as real, then  the  only
  $\mathfrak{sl}_2$-triple of  symmetric  type  is the  triple
  $(h,e,f)$ associated  with  the    real
    subalgebra  $\mathfrak{sl}_2\R$  of   the   Lie  algebra $\mathfrak{sl}_2\C$.
    The  associated canonical  decomposition~(\ref{main-decomposition})  is    given  by
\[ \mathfrak{sl}_2(\C)  = (\R h + \R e  + 0+ \R(i e )) +  (\R f   + (\R ih   + \R if)) \]
  \end{prop}

\Pf Denote by $J$ the $\ad{}$-invariant complex structure
on $\Lg$, and denote by $s$ the involutive automorphism that defines
$\Lg=\Lg_++\Lg_-$.
The complex structure $sJs$ is also
$\ad{}$-invariant, so simplicity of $\Lg$ implies
that $sJs=\pm J$. This means $J$ either commutes or anti-commutes
with~$s$. In the former case, $J\Lg_+=\Lg_+$ so $J\Ls$ would be a
simple ideal of $\Lg_+=\Ls+\Lz$. Since
$\Ls\cap J\Ls$ is a complex subalgebra of
$\Ls$ and $\Ls$ admits no non-trivial complex subalgebras,
we would have $J\Ls\subset\Lz$. However, this is a contradiction because
$\Lz$ centralizes $J\Ls$ and $J\Ls$ is not Abelian. 

We have seen that $J$ anti-commutes with $s$, so $J\Ls\subset\Lg_-$
and $\Ls\subset J\Lg_-\subset\Lg_+=\Ls+\Lz$. Recall that the adjoint
action of $\Ls$ admits no trivial components in $\Lg_-$, and
hence none on $J\Lg_-$. It follows that $\Ls=J\Lg_-$. 

The $\ad{}$-invariance of $J$ also implies that $J\Lz\subset\Lz$,
so $\Lz=J\Lz\subset\Lg_+\cap\Lg_-=0$. We have proved that
$\Lg=\Ls+J\Ls$ is the complexification of $\Ls\cong\mathfrak{sl}_2\R$. \EPf

\medskip

Note that the symmetric pairs $(\mathfrak{sl}_2\C,\mathfrak{sl}_2\R)$
and $(\mathfrak{sl}_2\R\oplus\mathfrak{sl}_2\R,\mathfrak{sl}_2\R)$
given in Propositions~\ref{simple} and~\ref{abs-simple}
have the same complexification and in a certain sense
are dual one to the other; the results of these propositions
are collected in Table~\ref{tab:1}.

\section{Symmetric contact manifolds of non-conical type}\label{sec:non-conical}

In this section, we prove Theorem~\ref{main4}. We first recall
some facts about symplectic symmetric spaces and refer to~\cite{Bi}
for more details on them. 

\subsection{Symplectic symmetric spaces}

A \emph{symplectic symmetric space} is a connected affine symmetric
manifold endowed $N$ with a symplectic structure
which is invariant under the geodesic symmetries. The transvection
group $G$ of a symplectic symmetric space
(i.e. the connected group generated by the the geodesic symmetries)
acts transitively on $N$. The symmetry at the basepoint normalizes
$G$ and induces on its Lie algebra $\Lg$ an involution~$s$
and hence a symmetric decomposition $\Lg=\Lk+\Lp$ into the $\pm1$-eigenspaces.
The symplectic structure at the basepoint yields an 
$\mathrm{ad}_{\mathfrak k}$-invariant non-degenerate $2$-form $\omega$ on $\Lp$,
and the triple $(\Lg,s,\omega)$ is called a \emph{symplectic
  symmetric Lie algebra}. Conversely, every symplectic symmetric
Lie algebra gives rise to a unique, up to isomorphism, simply-connected
symplectic symmetric space. 

\subsubsection{Symplectic symmetric spaces of a semisimple Lie group}
Let $N=G/K$ be a simply-connected
symplectic symmetric space of a connected semisimple
Lie group~$G$, where $K$ is connected.
The Whitehead lemmas imply that the invariant symplectic
form can be written as $d\theta$ for a unique $\theta\in\Lg^*$,
where $d$ is the Chevalley coboundary operator.  
This element is dual under the Killing form to an element $\xi$
in the center of $\Lk$. It follows that $N$ decomposes as a
product $G_1/K_1\times\cdots\times G_k/K_k$, where the $G_i$ are simple
Lie groups, $K_i=K_i'\cdot Z(K_i)$, $K_i'$ is a semisimple
Lie group, the center $Z(K_i)\cong T^1$, $\R$ or $\C^\times$,
and $\xi=\xi_1+\cdots+\xi_k$, where $\xi_i\in Z(\Lk_i)$ is non-trivial.  
Moreover, $N$ is an equivariant symplectic covering of the adjoint orbit of
$\xi$ in $\Lg$.

\subsubsection{Symplectic symmetric spaces of a complex simple Lie group
  and graded Lie algebras of depth one}\label{cx-sympl-symm}
Let $G$ be a complex simple Lie group with Lie algebra $\Lg$.
Fix a Cartan subalgebra and an ordering in the set of associated roots $R$.
Denote the system of simple roots by $\Pi$. 
Recall that the \emph{Dynkin mark} of a simple root $\alpha\in\Pi$
is its (necessarily positive, integral)
coefficient in the expression of the highest root as a linear combination
of simple roots.
Each $\alpha\in\Pi$ with Dynkin mark $1$ defines a gradation of depth one
\[ \Lg = \Lg^{-1}+\Lg^0+\Lg^1 =
\Lg(-R_1)+\Lg(R_0)+\Lg(R_1) \]
where $\xi=\frac12h_\alpha$ is the grading element, $h_\alpha$ is the coroot,
$R_0=\{\beta\in R\;|\;B(\beta,\alpha)=0\}$ and $R_1=R^+\setminus R^0$.
It turns out that every gradation of depth one of $\Lg$ arises in this
way, up to conjugation~(\cite[\S3.2.3]{CS}; see also~\cite{W}).
The associated symmetric decomposition
\[ \Lg=\Lk+\Lp=\Lg^{even}+\Lg^{odd}=\Lg^0+
(\Lg^{-1}+\Lg^1) \]
    defines a complex symmetric space $ G/ K$ with complex symplectic
    structure defined by the $\mathrm{ad}_{{\mathfrak k}}$-invariant
    $2$-form $\omega=d(B\circ\xi)$, whose kernel is $\Lk$.

    \subsubsection{Realifications of complex symplectic symmetric spaces}
    Let $N=G/K$ be a complex symplectic symmetric space with
    symplectic form~$\omega$ as in
    subsection~\ref{cx-sympl-symm}. Then $N$ can be considered as a
    real manifold with real invariant symplectic structure
    $\omega^r=a\Re\omega+b\Im\omega$, where $a$, $b\in\R$, $a^2+b^2>0$,
    and it becomes a (real) symplectic symmetric space.

    \subsubsection{Real forms of symplectic complex symmetric spaces} 
    Let $N=G/K$ be a real symplectic symmetric space of an
    absolutely simple Lie group $G$
    with symmetric decomposition
    \[ \Lg=\Lk+\Lp. \]
    The symplectic form corresponds to an $\mathrm{ad}_{\mathfrak k}$-invariant
    $2$-form $\omega$ on $\Lg$ with $\ker\omega=\Lk$. Moreover, $\omega$ is
    the differential of an $\mathrm{ad}_{\mathfrak k}$-invariant $1$-form
    $\theta$ which can be written as $\theta=B\circ\xi$ for some semisimple
    element $\xi$ in the center of $\Lk$ such that $\mathrm{ad}_\xi\Lp=\Lp$.
    Let $x$, $y\in\Lp^{\mathbb C}$ be eigenvectors with respective eigenvalues
    $\lambda$ and $\mu$. Then $[x,y]\in\Lk^{\mathbb C}$ is an eigenvector
    with eigenvalue $\lambda+\mu=0$. This shows that $\mathrm{ad}_\xi$
    has only two eigenvalues $\pm\lambda$, where $\lambda\neq0$ is either real
    or purely imaginary. In either case we get a gradation of depth $1$
    \[ \Lg^{\mathbb C} = \Lk^{\mathbb C}+\Lp^{\mathbb C} = (\Lg^\mathbb C)^0
    +((\Lg^{\mathbb C})^{-1}+(\Lg^{\mathbb C})^1), \]
    whose grading element $d$ is $\frac1{\lambda}\xi$ or $\frac i{\lambda}\xi$.
    Moreover, the Lie algebra is the real form of $\Lg^{\mathbb C}$ defined by an
    anti-involution $\sigma$ such that $\sigma(d)=\pm d$. We say that $\Lg$
    is \emph{consistent} with the depth $1$ gradation of $\Lg^{\mathbb C}$.

    Conversely, if we are given
    a real form $\Lg^\sigma$ of a complex simple Lie algebra
    $\Lg$, defined by an anti-involution $\sigma$ which is consistent with a
    depth one gradation of $\Lg$ in the sense it preserves the grading
    element~$d$
    up to sign, then the symmetric decomposition
    $\Lg=\Lk+\Lp$ restricts to a symmetric decomposition
    $\Lg^\sigma=\Lk^\sigma+\Lp^\sigma$
    and the element $\xi\in\Lg^\sigma$ equal to $d$ or $i\cdot d$ defines
    an $\mathrm{ad}_{{\mathfrak k}^\sigma}$-invariant $2$-form $\omega=d(B\circ\xi)$
    on $\Lg^\sigma$, with kernel $\Lk^\sigma$. The corresponding symmetric
    manifold $N=G^\sigma/K^\sigma$ is a symplectic symmetric space with symplectic
    form induced by $\omega$. 

    \subsection{Symmetric contact spaces of non-conical type as
      $1$-dimensional bundles over symplectic symmetric spaces}

    Let $(M=G/H,\D)$ be a homogeneous contact manifold of a
    connected semisimple Lie group $G$ which is the total
    space of a $1$-dimensional bundle over $N=G/K$
    as in subsection~\ref{non-conical}. Assume $(M,\D)$ is a
    symmetric contact space, so there exists an involutive
    automorphism $s$ of $\Lg$ that preserves $\Lh$ and induces
    $-1$ on $\D$. 
    
\begin{prop}\label{sym-non-con}
There exists a symmetric decomposition
\begin{equation*}\label{sympl-inv}
 \Lg = \Lk + \Lp 
\end{equation*}
under~$s$, where
\[ \Lk = \Lh + \R\eta \]
and~$\eta$ is a nonzero element in the center of $\Lk$.
Moreover, $\mathcal D$ is the $G$-invariant extension of $\Lp$ and
there exists an $\ad{\mathfrak k}$-invariant symplectic
structure on $\Lp$ given by $d\theta$, where $\theta\in\Lg^*$,
so that $(\Lg,s,d\theta)$  
is a symplectic symmetric Lie algebra. Finally, $\eta$
generates a central closed subgroup $C$ of~$K$ and
$K=C\times H$ (direct product).
\end{prop}

\Pf Recall that
\[  \Lg = \Lk  + \Lp,\quad \Lk = \Lh + \R\eta, \]
where $\Lk$ is the centralizer of $\theta\in\Lg^*$, 
$\Lh=\ker\theta\cap\Lk$ is a codimension one ideal of $\Lk$ and
$\theta(\eta)=1$ (cf.~(\ref{decomp-non-con})).
Since $s$ must preserve the contact structure,
$s^*\theta=c\,\theta$ for some $c\neq0$, and
the involutivity forces $c=\pm1$. Now $s$ is a
semisimple automorphism of $\Lg$ that
preserves $\ker\theta$ and $\Lk=\ker d\theta$,
so we can choose the subspace $\Lp$ 
to be $s$-invariant. 

We have $\ker\theta=\Lh+\Lp$
so $\Lp$ induces the distribution $\mathcal D$ and $s|_{\mathfrak p}=-1$.
Let $\xi\in\Lg$ be the element
dual to $\theta$ under the Killing form $B$.
Then
\begin{equation}\label{xi-w}
 0=\theta(\Lh+\Lp)=B(\xi,\Lh+\Lp). 
\end{equation}
Plainly, $\xi\in Z_{\mathfrak g}(\xi)=Z_{\mathfrak g}(\theta)=\Lk$, so
$\xi$ lies in the center of $\Lk$. 

We already know that $s\xi=\pm\xi$. Non-degeneracy of $d\theta$ on $\Lp$
gives $x$, $y\in\Lp$ such that $d\theta(x,y)=-1$ and then
\[ s[x,y]=s(\eta+[x,y]_{\mathfrak p}+[x,y]_{\mathfrak h})=s\eta-[x,y]_{\mathfrak p}+(s[x,y])_{\mathfrak h} \]
and
\[ [sx,sy]=[-x,-y]=\eta+[x,y]_{\mathfrak p}+[x,y]_{\mathfrak h} \]
implying that $s\eta=\eta$. Hence $s^*\theta=\theta$ and $s\xi=\xi$.
It also follows that 
\begin{equation}\label{eta-w}
B(\eta,\Lp)=0.
\end{equation}

Write $\Lg=\Lg_++\Lg_-$ and $\Lh=\Lh_++\Lh_-$ under~$s$.
Note that $\Lg_+=\Lh_++\R\eta$ is a reductive
subalgebra of $\Lg$ (this follows for instance from~\cite[Lemma~20.5.12]{T-Y})
and $\Lg_-=\Lh_-+\Lp$, $\ker\theta=\Lh_++\Lg_-$. 
Since $[\Lk,\Lg]\subset\Lh$, 
we have $[\Lh_-,\Lg_+]\subset\Lh_-$. Therefore
we can change $\Lp$, if necessary, so that $\Lp$ is in addition
$\mathrm{ad}_{\mathfrak g_+}$-invariant.

We next claim $\ad\xi:\Lp\to \Lp$ is an isomorphism. In fact, if $\ad\xi x=0$ for
some $x\in \Lp$, then $d\theta(x,\Lp)=B(\ad\xi x, \Lp)=0$ implying $x=0$ by
nondegeneracy of $d\theta$ on $\Lp$.

It is clear that $B$ is nondegenerate on $\Lg_+\subset\Lk$.
Next, note that $[\Lh_-,\Lp]\subset\Lg_+$ and
$B([\Lh_-,\Lp],\Lg_+)=B(\Lp,[\Lh_-,\Lg_+])\subset B(\Lp,\Lh_-)=B(\ad\xi \Lp,\Lh_-)=B(\Lp,\ad\xi\Lh_-)=0$, since $\xi$ centralizes~$\Lh$.
We have proved that $[\Lh_-,\Lp]=0$. It follows that
there is a symmetric decomposition
\[ \Lg=(\Lg_++\Lp)+\Lh_-. \]
This argument also shows that 
\begin{equation}\label{h-w}
B(\Lh,\Lp)=B(\Lh,\mathrm{ad}_\xi \Lp)=B(\mathrm{ad}_\xi\Lh,\Lp)=0.
\end{equation}
Denote by $\Ln_3$ the kernel of the adjoint action of $\Lg_++\Lp$ on $\Lh_-$;
this is an ideal of $\Lg$ that contains $\R\xi+\Lp$. Using~(\ref{xi-w}),
(\ref{eta-w}) and~(\ref{h-w}), we deduce that
\[ \Ln_3^\perp\subset \Lp^\perp\cap\xi^\perp=(\Lh+\R\eta)\cap(\Lh+\Lp)=\Lh, \]
and hence $\Ln_3^\perp=0$, owing to the effectiveness of the
presentation $M=G/H$.
Finally $\Ln_3=\Lg$ and thus $\Lh_-=0$, so that we arrive at
\[ \Lg=\Lk + \Lp, \]
symplectic involutive Lie algebra, where $\Lk=\Lh+\R\eta=\Lg_+$ 
and the symplectic structure on $\Lp=\Lg_-$ is induced by $d\theta(x,y)=-B(\ad\xi x,y)$ for $x$, $y\in \Lp$. Since $\Lk=\Lg^+$ is a reductive subalgebra of $\Lg$,
$\Lk=[\Lk,\Lk]\oplus Z(\Lk)$, where $Z(\Lk)$ denotes the
center of~$\Lk$, $[\Lk,\Lk]\subset\Lh$, and 
$\eta$ can be chosen in the center of $\Lk$ and to generate 
a closed subgroup of $G$.  
\EPf

\medskip

The next example shows that $\xi$ and $\eta$
in Proposition~\ref{sym-non-con} do not have to coincide. 

\begin{ex}\label{ex}
  Consider the complex Lie algebra $\mathfrak{sl}_2(\C)$
  with its standard basis $(h,e,f)$.
  Its Killing form $B^c$ satisfies $B^c(h, h) =8$
  and $B^c(e, f ) = 4$ with other values zero.
  Let $\Lg$ denote the realification of $\mathfrak{sl}_2(\C)$.
  The Killing form of $\Lg$ is $B = 2\Re B^c$.
  We obtain a symmetric contact space of non-conical type by
  choosing $\xi=\lambda h\in\Lg$, where $0\neq\lambda\in\mathbb C$.
  In fact $\Lg = \Lk + \Lp$
  is a symmetric decomposition, where $\Lk = \C h$, $\Lp = \C e + \C f$
  and $\Lh = \R(\mu h)$, where $\Re(\lambda\mu)=0$.
  Note that $\xi$ is $B$-isotropic if and only
  $\Re(\lambda^2)=0$; in that case, $\eta\neq\xi$. 
  \end{ex}

\textit{Proof of Theorem~\ref{main4}.}
Let $N=\Ad G\xi$ be a symplectic symmetric space as in the statement. 
Then $N=G/K$, where $K$ is connected, and there exists a symmetric
decomposition $\Lg=\Lk+\Lp$ under an involution $s$. 
Moreover the symplectic structure on $N$ is induced from
the $\mathrm{Ad}_K$-invariant form given by the restriction
to~$\Lp$ of
$\omega=d\theta$, where $\theta\in\Lg^*$
is dual to~$\xi$ under the Killing form $B$ of $\Lg$.
Let $\Lh=\ker\theta\cap\Lk$ ($\Lh$ is also the $B$-orthogonal of $\xi$ in $\Lk$)
and consider the associated connected subgroup $H$ of $K$.
It turns out $H$ is a codimension one normal subgroup of $K$. 
By assumption, $H$ is closed in $K$, so
$M=G/H$ is the total space of a principal 
$K/H$-bundle over $N$. The $1$-form $\theta$ defines an invariant
contact structure on $M$ which is the invariant extension of
the hyperplane $\Lp\subset T_oM$. The involution~$s$ preserves $\Lh$
and $\theta$, thus it induces a contactomorphism of $M$. 
Moreover, it induces $-1$ on $\mathcal D=(\ker\theta)/\Lh$ and
hence $M$ is a symmetric contact space.

The converse follows from Proposition~\ref{sym-non-con}. \EPf

\begin{rem}
  In cases $\dim Z(\Lk)=1$ or $\dim_{\mathbb C} Z(\Lk)=1$
  (i.e.~$\Lg$ is absolutely simple or complex simple;
  such cases are exactly those listed in Tables~\ref{tab:5},
  \ref{tab:6}, \ref{tab:7} and~\ref{tab:8}),
  $H$ is closed for any choice of $\xi\in Z(\Lk)$;
  in the first case, in addition,
  $\eta$ can be taken to be a multiple of $\xi$.
  In general, $H$ is closed if and only if $H\cap Z(K)^0$ is closed;
  a sufficient condition is that there exists a non-compact (closed)
  one-parameter subgroup of $K$ not contained in $H$.   
\end{rem}


\vfill

\end{document}